\documentclass[11pt]{amsart}
\usepackage{graphicx}
\usepackage[all]{xy}
\usepackage{latexsym}
\usepackage{color}
\usepackage{xcolor}
\usepackage{amsmath,amssymb,amsthm,mathtools}
\xyoption{arc}

\usepackage{tikz, tikz-cd}
 \tikzstyle{int}=[circle, draw,fill=black,outer sep=0,minimum size=3pt, inner 
  \tikzstyle{ext}=[circle, draw=black,outer sep=0,inner sep=1pt]

\usepackage{fullpage}

  \newcommand{\SV}{{\Sigma,V}}

   \newcommand{\FR}{\mathsf{FR}}

\newcommand{\id}{\mathrm{id}}

 \newcommand{\lon}{\longrightarrow}
 \newcommand{\bu}{\bullet}
 \newcommand{\ad}{{\mathrm a\mathrm d}}
 \newcommand{\rar}{\rightarrow}
 
 \newcommand{\End}{\mathsf{End}}
\newcommand{\p}{{\partial}}
  \newcommand{\ot}{\otimes}
 \newcommand{\wh}{\widehat}

\newcommand{\Der}{\mathrm{Der}}

 \newcommand{\Hom}{{\mathrm{Hom}}}
 \newcommand{\Mor}{{\mathrm{Mor}\,}} 
 
 \newcommand{\sip}{\smallskip}
 \newcommand{\bip}{\bigskip}
 \newcommand{\mip}{\vspace{2.5mm}}
 
\newcommand{\Lie}{\mathcal{L} \mathit{ie}}
\newcommand{\BV}{{\mathcal B}{\mathcal V}}

 \newcommand{\Z}{{\mathbb Z}}

 \newcommand{\R}{{\mathbb R}}
 \newcommand{\N}{{\mathbb N}}
 \newcommand{\K}{{\mathbb K}}

\newcommand{\sS}{{\mathsf S}}

 \newcommand{\Beq}{\begin{equation}}
 \newcommand{\Eeq}{\end{equation}}
 \newcommand{\Beqr}{\begin{eqnarray}}
 \newcommand{\Eeqr}{\end{eqnarray}}
 \newcommand{\Beqrn}{\begin{eqnarray*}}
 \newcommand{\Eeqrn}{\end{eqnarray*}}
 \newcommand{\Ba}{\begin{array}}
 \newcommand{\Ea}{\end{array}}
 \newcommand{\Bi}{\begin{itemize}}
 \newcommand{\Ei}{\end{itemize}}
 \newcommand{\Bc}{\begin{center}}
 \newcommand{\Ec}{\end{center}}

 \newcommand{\fe}{{\mathfrak e}}

 \newcommand{\fg}{{\mathfrak g}}
 \newcommand{\fh}{{\mathfrak h}}

 \newcommand{\fG}{{\mathfrak G}}

\newcommand{\fL}{{\mathfrak L}}
\newcommand{\fM}{{\mathfrak M}}

 \newcommand{\f}{{\mathcal O}}
 \newcommand{\cA}{{\mathcal A}}
 \newcommand{\cB}{{\mathcal B}}

 \newcommand{\cF}{{\mathcal F}}

 \newcommand{\caL}{{\mathcal L}}
 \newcommand{\cM}{{\mathcal M}}

 \newcommand{\cT}{{\mathcal T}}

 \newcommand{\cV}{{\mathcal V}}

\newcommand{\Sym}{\mathrm{Sym}}
\newcommand{\Fun}{C^\infty}

\newcommand{\Hol}{\mathrm{Hol}}

 \newcommand{\al}{\alpha}
 \newcommand{\be}{\beta}
 \newcommand{\ga}{\gamma}
 
 \newcommand{\Ga}{\Gamma}

 \newcommand{\la}{\lambda}
 \newcommand{\om}{\omega}

\theoremstyle{plain}
\swapnumbers

\newtheorem{prop-def}[theorem]{Proposition-definition}

\newtheorem{f-theorem}{Formality Theorem}[section]
\newtheorem{main-theorem}{Main~Theorem}[section]
\newtheorem{section-theorem}{Theorem}[section]

\theoremstyle{definition}

\begin{document}

\title{Homotopy Poisson and BV$_\infty$ structures \\
on derived moduli spaces of flat connections}
\author{Sergei~A.\ Merkulov}

\begin{abstract}
We show that quasi-Poisson and quasi-BV structures on (super)manifolds $M$ admitting an action of a Lie (super)group $G$ can be identified with a special class of homotopy  Poisson and respectively $BV_\infty$ structures  on the associated differential (super)manifolds of the form $\cM:=M\times \Pi\fg$, $\fg$ being the Lie (super)algebra of $G$. The construction is based on the proven existence of an essentially unique family of higher Jacobiators for any Lie (super)algebra $\fg$ with an invariant (odd or even) scalar product.
 We apply these arguments to the moduli spaces $M$ of flat $G$-connections on a surface with boundaries using Fock-Rosly construction.


\end{abstract}

 \maketitle
 \markboth{}{}

{\Large
\section{\bf Introduction} 
}

\mip

\subsection{From quasi-Poisson to dg Poisson} Let $M$ be a smooth manifold $M$ equipped with an action of a Lie group $G$ whose Lie algebra $\fg$ has an invariant scalar product.
The first main result of this paper is an observation that any quasi-Poisson structure $\pi$  on  $M$  \cite{AKSM} can be equivalently understood as a $Poisson_\infty$ structure on a  graded (super)manifold
$$
\cM:= M \times \fg[1],
$$
whose odd sector equals precisely the degree shifted Lie algebra. This graded manifold
comes canonically equipped with a degree +1 vector field $Q$ such that
$Q^2=0$, and the associated cohomology group $H^0 (\Fun(\cM),Q)$ is equal to the algebra of the $G$-invariant functions $\Fun(M)^G$ on the original manifold $M$.

\mip

The aforementioned  $Poisson_\infty$-structure  on $\cM$ is of a special type: it is best understood as an ordinary graded Poisson structure (that is, an MC element $\mu_\hbar$ in the graded Lie algebra of polyvector fields on $\cM$) of the form 
$$
\mu_\hbar=Q +  \hbar\pi + \sum_{n\geq 1} \hbar^{2n}\mu_{2n+1}\in T\cM \oplus \wedge^2 T\cM[[\hbar]] \oplus \prod_{n\geq 1}\wedge^{2n+1} T\cM[[\hbar]], \ \ \ 
[\mu_\hbar, \mu_\hbar]_{\mathrm{Schouten}}=0.
$$
Moreover, our construction of higher homotopies $\{\mu_{2n+1}\}_{n\geq 1}$ is universal, it depends on the Lie algebra $\fg$, its representation 
$$
\iota: \fg \lon TM, 
$$
 and its scalar product; it does not depend on the choice of a particular quasi-Poisson structure $\pi$.
  
  \sip
  
 Plainly speaking, the idea to get ``dg Poisson out of quasi-Poisson''  works fine due to the existence of a universal family of  so-called {\it higher Jacobiators},
$$
\left\{\nu_{2n+1} \in
 \left(\wedge^2\fg
 \ot \odot^{2n-1}\fg \right)^{\fg}\right\}_{n\geq 1}
$$
with $\nu_3\in (\wedge^2 \fg\ot \fg)^\fg$ being uniquely determined by the standard Cartan trivector $\phi\in (\wedge^3\fg)^\fg$ which is used in the definition \cite{AKSM} of quasi-Poisson  structures on $\fg$-manifolds. Their representations in  $\wedge^{2n+1} T\cM$ induced by the given representation $\iota$ 
are precisely the higher components 
$$
\mu_{2n+1}=\iota(\nu_{2n+1}), \ \ \forall n\geq 1,
$$
of the graded Poisson structure on $\cM$. Moreover, the family of higher Jacobiators for any Lie algebra $\fg$ with an invariant scalar product is defined essentially uniquely (see \S {\ref{2: subsec on essential uniqueness}} for a precise statement). As an illustration, the list of {\it all}\, higher Jacobiators in the case $\fg=\mathfrak{sl}_2(\K)$ is given in \S {\ref{2: example sl_2}}.

\sip

The graded Poisson structure $\mu_\hbar$ can be deformation quantized via the classical Kontsevich formality  map \cite{Ko} making the graded commutative algebra $\Fun(\cM)[[\hbar]]$ into a (curved, in general) $A_\infty$-algebra. Note that this quantization idea applies to a generic quasi-Poisson manifold, not to a special class of quasi-Poisson manifolds whose deformation quantizations have been studied earlier in \cite{EE,LS,H} via  different approaches.

\subsection{From quasi-BV to BV$_{com}^\infty$} Let $M$ be a supermanifold  equipped with an action of a Lie supergroup $G$ whose Lie algebra $\fg$ has an invariant odd scalar product.
The second main result of this paper says that any  quasi-BV \cite{ANPS}  structure $\Delta$ on $M$ can be equivalently understood as $BV_\infty^{com}$ structure $\Delta_\hbar$,
$$
\Delta_\hbar=Q +\hbar \Delta + \sum_{n\geq 1} \hbar^{2n}\Delta_{2n+1}: \Fun(\cM) \rar \Fun(\cM)[[\hbar]], 
\ \ \  \ \Delta_\hbar^2=0,
$$
on the associated dg supermanifold
$$
\cM:= M \times \Pi \fg.
$$
This supermanifold is equipped with an odd vector field $Q$ such that
$Q^2=0$ and  $\Fun(M)^G$  sits in  $H^{even} (\Fun(\cM),Q)$. The odd operators
 $\Delta_{2n+1}$ on $\Fun(\cM)$ are of order $2n+1$.
Our construction of higher terms $\{\Delta_{2n+1}\}_{n\geq 1}$ of the $BV_\infty^{com}$ structure is universal, it depends only on the Lie superalgebra $\fg$, its representation in $TM$ and its odd scalar product; it does not depend on the choice of a particular pre-BV structure $\Delta$.

\sip

An important class of quasi-Poisson and quasi-BV manifolds comes from the theory of framed moduli spaces $M$ of flat $G$-connections on surfaces $\Sigma$ of any genus $g$ with $m\geq 1$ boundaries \cite{ANPS, FR, LS, Ma, MT, Mo, N}. They are best understood and explicitly described  via the famous Fock-Rosly construction \cite{FR} based on a choice of a ribbon skeleton $\Ga$ of $\Sigma$. We follow the same approach and construct the above mentioned extensions $\mu$ and $\Delta_\hbar$ on the associated graded manifolds $\cM$ of ``derived'' flat connection  via a choice of a skeleton $\Ga\subset \Sigma$, and then observe that the result does not depend on that choice, i.e.\ the resulting $Poisson_\infty$ and $BV_\infty^{com}$ structures are intrinsic.

\sip

\subsection{Structure of the paper} In \S 2 we show that every quasi-Poisson structure on a $G$-manifold $M$ can be understood as a graded Poisson structure on a dg manifold $\cM=M\times \fg[1]$, $\fg$ being the Lie algebra of the given Lie group $G$ (see Main Theorem {\ref{2: Main Theorem on mu_hbar}}). In \S 3 we apply this result to the framed moduli spaces of flat $G$-connections on a surface $\Sigma$ with boundary. In \S 4 we show that every quasi-BV structure on a $G$-supermanifold $M$ can be understood as a $BV_\infty^{com}$ structure on a differential supermanifold $\cM=M\times \Pi\fg$ (see Main Theorem {\ref{4: Theorem on Delta_hbar}}). Following \cite{ANPS} this result applies to the framed moduli space of flat $G$-connections on $\Sigma$ such that the Lie superalgebra $\fg$ of the given Lie supergroup $G$ has an odd invariant scalar product.

\subsection{Some notation} We work over a field $\K$ of characteristic zero.
The cardinality of a finite set $I$ is denoted by $\# I$; its linear span over a field $\K$ is denoted by $\K[I]$. If $V=\oplus_{i\in \Z} V^i$ is a $\Z$-graded vector space and $n\in \Z$, then
$V[n]$ stands for the graded vector space with $V[n]^i:=V^{i+n}$; for $v\in V^i$ we set $|v|:=i$.

\sip
{\bf Acknowledgements}. Partial support from the Basic Research Program of HSE University is gratefully acknowledged (HSE-BR-2025-84). 

\sip

{\bf Declaration of AI usage}.
The author used ChatGPT (OpenAI, model GPT-5.6 Sol Light) to assist with proof-reading, checking arguments for possible
errors, and carrying out some routine but lengthy calculations in
Example~{\ref{2: example sl_2}}. These calculations were subsequently
checked by the author. Apart from this explicitly stated assistance,
all mathematical ideas, results, and proofs were developed and verified
by the author, who assumes full responsibility for the entire content
of the manuscript.

\bip


{\Large
\section{\bf Quasi-Poisson structure as a dg Poisson structure}
}

\mip

\subsection{Notation for a Lie algebra with an invariant scalar product}\label{2: subsec about G and fg}
The main results of this paper about $Poisson_\infty$ and $\BV^{com}_\infty$ structures associated with quasi-Poisson and quasi-BV (super)manifolds
are proven by using explicit coordinate representations of the aforementioned homotopy structures. 
Hence we start with a brief summary of
the conventions and notation for all the structure constants involved in the calculations of $Poisson_\infty$-structures. A similar summary for the case of Lie superalgebras with an odd scalar product --- which are used in the study of quasi-BV supermanifolds --- will be given below. 

\sip

Let $G$ be a connected Lie group with  Lie algebra $(\fg, [\ ,\ ])$.
The $Poisson_\infty$ structure on the ``derived" moduli space of flat $G$-connections on a closed surface $\Sigma$ with boundaries is described in this paper in terms of explicit polydifferential operators. These operators are best described in coordinates associated with some basis $\{e_i\}_{i\in I}$ of $\fg$.
Moreover, the proofs of the computability relations between these polydifferential operators  are given by straightforward calculations in these coordinates.

\sip

The structure constants of the Lie bracket on $\fg$ are given in the basis  $\{e_a\}_{a\in I}$  by
\Beq\label{2: Lie bracket in fg}
[e_a, e_b]=\sum_{c} C_{ab}^c e_c.
\Eeq
 We assume that $\fg$ has a symmetric non-degenerate invariant tensor 
\Beq\label{2: invariant tensor f}
 t=\sum_{a,b} t^{ab} e_a\otimes e_b \in \odot^2\fg
\Eeq
The invariance of $t$ means
$$
 [t, X\otimes \id +\id \otimes X]=0
 \ \ \ \forall X\in\fg
$$
or, in terms of the structure constants,
$$
\sum_{k}\left(t^{ak}C_{kc}^b + t^{bk}C_{kc}^a\right)=0.
$$
The non-degeneracy of $t$ means that the linear map
$$
\Ba{rccc}
t^\#: & \fg^* & \lon & g\\
   & X'=\sum_a x'_a e^a & \lon & 
   \sum_{a,b} (t^{ab} x'_a) e_b 
\Ea
$$
is an isomorphism. Hence one has a non-degenerate scalar product
$$
\Ba{rccc}
\langle\ ,\ \rangle: & \odot^2\fg & \lon & \K\\
   & (X=\sum_a x^a e_a, Y=\sum_b y^b e_b) & \lon & 
   \langle X,Y \rangle:= \sum_{a,b} t_{ab} x^ay^b 
\Ea
$$
 on $\fg$, where $t_{ab}$ is the matrix inverse to $t^{ab}$. The associated Cartan trivector $\phi^*\in \wedge^3\fg^*$ 
 is a map
 $$
\Ba{rccc}
\phi^*: & \wedge^3\fg & \lon & \K\\
   & (X, Y,Z) & \lon & 
   \phi^*(X,Y,Z):= \langle X, [Y,Z]\rangle 
\Ea
$$
 The isomorphism $t^\#: \fg^* \simeq \fg$ gives us an associated dual element
 $\phi\in \wedge^3 \fg$ given in coordinates by
\Beq\label{2: cartan trivector}
\phi
=\sum_{a,b,c}\phi^{abc}
e_a\wedge e_b\wedge e_c \in \wedge^3 \fg,
\Eeq
where
\Beq\label{2: phi_abc}
\phi^{abc}=\sum_{k,l} t^{ak}t^{bl}C_{kl}^c.
\Eeq
The invariance of $t$ implies skew-symmetry of the structure constants $\phi^{abc}=\phi^{[abc]}$.
The Cartan trivector is itself $\fg$-invariant in the sense that
the equation
\Beq\label{2: inv of Cartan trivector}
\sum_k \left( \phi^{kbc}C_{kl}^a +  \phi^{kab}C_{kl}^c  + \phi^{kca}C_{kl}^b\right)=0
\Eeq
holds for any $a,b,c\in I$.

\subsection{Quasi-Poisson manifolds}
Let  $G$ be a connected  Lie group with Lie algebra $\fg$ equipped with an invariant scalar product as in the previous subsection. Assume $G$ acts smoothly on a manifold $M$. Differentiating the action gives a Lie algebra homomorphism
\Beq\label{2: iota}
\iota:\fg\longrightarrow\Gamma(TM).
\Eeq
This map extends naturally to exterior powers,
$$
\Ba{rccc}
\iota: & \wedge^k\fg &\lon &
\Ga(\wedge^kTM)\\
X_1\wedge\cdots\wedge X_k &\lon & \iota(X_1)\wedge\cdots\wedge\iota(X_k)
\Ea
$$
In particular, the Cartan trivector 
(\ref{2: cartan trivector}) gives a trivector field
$$
\iota(\phi)\in (\wedge^3TM).
$$
on the manifold $M$.

\sip

A polyvector field $\nu\in \Ga(\wedge^\bu TM)$ is called {\it $G$-invariant}\, if 
$$
[\iota(X), \nu]_\sS=0 \ \ \ \forall X\in \fg,
$$
where $[\ ,\ ]$ is the standard Schouten bracket of polyvector 
fields. In particular, $G$-invariant polyvector fields $\nu$ form a Lie subalgebra $(\wedge^\bu TM)^G$
of the Lie algebra $(\wedge^\bu TM, [\ , \ ]_\sS)$, and have the property that
$$
(L_g)_*(\nu)|_{g\cdot x}=\wedge^\bu(dL_g)_x (\nu|_x)
$$
for any $x\in M$ and the diffeomorphism $L_g:M\rar M$ sending $x$ to $g\cdot x$, $\forall \ g\in G$.

\sip

A {\it quasi-Poisson structure}\, on a $G$-manifold $M$ is, by definition \cite{AKSM}, a $G$-invariant bivector $\pi\in (\wedge^2 TM)^G$ satisfying the condition
\Beq\label{2: quasi-Poisson def}
\frac{1}{2}[\pi,\pi]_{\sS}=\iota(\phi).
\Eeq

Framed moduli spaces of trivial $G$-bundles over a surface (with possibly boundary) gives us an important class of examples of quasi-Poisson structures. We study below a ``derived" version of this class in much detail.

\subsection{The Chevalley--Eilenberg model of the algebra of invariant functions}\label{2: CE complex of G-manifolds}
The morphism of Lie algebras (\ref{2: iota}) makes the commutative algebra $\Fun(M)$ of smooth functions on the manifold $M$ into an $\fg$-module. Consider the associated  Chevalley--Eilenberg 
cochain complex
$$
C^\bullet(\fg,\Fun(M)):=\Hom_\R(\wedge^\bu g,  C^\infty(M))\simeq 
C^\infty(M)\otimes \Lambda^\bullet \fg^*
$$
equipped with the differential $Q$ 
$$
C^\infty(M) \stackrel{Q}{\lon} \Hom_\R(g,  C^\infty(M)) 
\stackrel{Q}{\lon} \Hom_\R(\wedge^2 g,  C^\infty(M))  \stackrel{Q}{\lon}\ldots
$$
which, by definition, acts  on a $p$-th cochain
$$
\Ba{rccc}
\om: & \wedge^p\fg &\lon &  C^\infty(M)\\
     & (X_1,\ldots, X_p) & \lon & \om (x_1,\ldots, X_p)
\Ea
$$
as follows
\begin{align*}
(Q \om)(X_0,\ldots,X_p)
=&
\sum_{i=0}^p (-1)^i
\rho(X_i)
\omega(X_0,\ldots,\widehat X_i,\ldots,X_p)
\\
&+
\sum_{0\leq i<j\leq p}
(-1)^{i+j}
\omega([X_i,X_j],
X_0,\ldots,\widehat X_i,\ldots,\widehat X_j,\ldots,X_p).
\end{align*}
The zero-th cohomology of this complex is precisely the algebra of $G$-invariant functions
$$
H^0(\mathcal \cA,\mu_1)
= \Fun(M)^G,
$$
the ones, one which $\fg$ acts trivially,
$$
\Fun(M)^G:=\{f\in \Fun(M) |\ \iota(X)f=0 \ \ \forall X\in \fg\}.
$$
The a quasi-Poisson structure $\pi$ on $M$ restricted to the subalgebra 
$\Fun(M)^G$ becomes a Poisson structure as for any $f,g,h \in \Fun(M)^G$
one has
$$
\iota(\phi)(df, dg,dh) = \sum_{a,b,c}\phi^{abc}
\bigl( \iota(e_a)f \bigr) \bigl( \iota(e_b)g\bigr)\bigl( \iota(e_c)h\bigr)=0.
$$
Hence the smooth part (if any) of the space of orbits $M/G$ is a Poisson manifold, but it is often a  bad idea to work with this quotient. We prefer in this paper to work with the full complex $C^\bullet(\fg,\Fun(M))$ as a kind of ``derived'' model for $M/G$; it can be understood super-geometrically as follows.

\subsection{$C^\bullet(\fg,\Fun(M))$ as a dg manifold}\label{2: CE as cM} Assume $\dim M=m$ and $\dim G=n$.
To simplify our calculations, it is useful to understand the graded commutative algebra $C^\infty(M)\otimes \Lambda^\bullet \fg^*$ as an algebra of smooth functions on the following $(m|n)$ dimensional {\it graded}\, (super)manifold
$$
\cM:= M \times \fg[1]
$$
where elements of the factor $\fg[1]$ are assumed to sit in degree $-1$.
The dual basis $\{e^a\}$ of $\fg^*$ induces odd (of degree +1) coordinates on $\cM$ which we denote from now on by $\{\eta^a\}_{a\in [n]}$,
$$
\eta^a\eta^b=-\eta^b\eta^a.
$$
If we denote (arbitrary) local coordinates on the even part $M$ by $x=\{x^\al\}_{\al\in [m]}$, then every smooth function on the supermanifold\footnote{Since we use $\Z$-grading rather than $\Z_2$-grading, we call $\cM$ a {\it graded}\, manifold from now on rather than a supermanifold.} $\cM$ is given by a finite sum
$$
f(x,\eta)=\sum_{k\geq 0} \sum_{a_1,..., a_k} \frac{1}{k!} f_{a_1...a_k}(x)\eta^{a_1}\cdots \eta^{a_k}.
$$
The odd left partial derivative $\p/\p \eta^a$ acts on such functions in the way. In this notation the differential $Q$ is given by the following odd vector field 
\Beq\label{2: Q}
Q= \sum_{a=1}^n \eta^a \iota(e_a) - \frac{1}{2} \sum_{a,b,c=1}^n \eta^a \eta^b C_{ab}^c \frac{\p}{\p \eta^c}
\Eeq
on the graded manifold $\cM$. In this form it is almost obvious that $\mu_1$ is a derivation of the graded commutative algebra $C^\infty(\cM)$, and satisfies the condition 
$$
[Q,Q]_\sS=2 Q^2=0.
$$
Hence the pair $(\cM,Q)$  is a dg manifold. Here $\frac{\p}{\p \eta^c}$ is the left partial derivative uniquely specified by the condition
$$
\frac{\p}{\p \eta^c}(\eta^a)=\delta_c^a.
$$

\subsection{From a quasi-Poisson structure on $M$ to a graded Poisson structure on $\cM$}\label{2: subsec from pi to mu}
Let $\wedge^\bu T\cM$ stand for Gerstenhaber algebra of polyvector fields on the graded manifold $\cM$ equipped with the standard Schouten bracket\footnote{We use the unique grading of $\wedge^\bu T\cM$ in which the standard Schouten Lie bracket $[\ ,\ ]_\sS$ has precisely degree $-1$, while the standard wedge  multiplication in $\wedge^\bu T\cM$ has degree zero. In this grading quasi-Poisson structures $\pi\in \wedge^2 T\cM\subset \wedge^2 T\cM$ have degree 2 exactly.} , and let $\pi\in (\wedge^2TM)^G$ be any quasi-poisson structure on the original manifold $M$. Let $\hbar$ be a degree zero formal parameter. Let us associate to $\pi$ the following polyvector field on $\cM$ (of total degree +1)
$$
\mu_\hbar^{< 4}:=\mu_1 +\hbar\mu_2 +\hbar^2\mu_3 \in T\cM \oplus \wedge^2T\cM[[\hbar]] \oplus \wedge^3 T\cM[[\hbar]]
$$
where, in an arbitrary local coordinate system $(x^\al, \eta^a)$ on $\cM$,
the vector field $\mu_1$ is given by (\ref{2: mu_1}) and
\begin{align}
\mu_1&:=Q \equiv \sum_{a=1}^n \eta^a \iota(e_a) - \frac{1}{2} \sum_{a,b,c=1}^n \eta^a \eta^b C_{ab}^c \frac{\p}{\p \eta^c}, \label{2: mu_1} \\
\mu_2&:=\pi= \sum_{\al,\be\in [m]} \pi^{\al\be}(x) \frac{\p}{\p x^\al}\wedge 
\frac{\p}{\p x^\be}, \label{2: mu_2}
\\
\mu_3 &:= -\sum_{a,b,c\in [n]} \phi^{abc}
\iota(e_a) \wedge \iota(e_b) \wedge \frac{\p}{\p \eta^c}.
\label{2: mu_3}
\end{align}

\subsubsection{\bf Lemma}\label{2: Lemma on 3-term mu Poisson} {\it The polyvector field $\mu_\hbar^{<4}\in \wedge^\bu T\cM[[\hbar]]$ satisfies}
$$
[\mu_\hbar^{<4},\mu_\hbar^{<4}]_\sS=0 \mod O(\hbar^4)
$$

\begin{proof}
The theorem is proven once we establish the following set of three identities,
$$
 [\mu_1, \mu_2]_\sS\equiv 0, \ \ [\mu_1,\mu_3]_\sS+ \frac{1}{2}[\mu_2,\mu_2]_\sS\equiv 0, \ \ 
[\mu_2,\mu_3]_\sS\equiv 0.
$$
The identities  $[\mu_1, \mu_2]_\sS\equiv 0$ and $[\mu_3,\mu_2]\equiv 0$ are immediate because 
$\pi$ is $G$-invariant, i.e.\ $[\iota(e_a),\pi]_\sS=0$ for any $a\in [n]$, and $\pi$ does not depend on odd coordinates. 

\sip

Next we have
\Beqrn
-[\mu_1,\mu_3]_\sS&=&\left[\sum_{l=1}^n \eta^l \iota(e_l) - \frac{1}{2} \sum_{p,q,k=1}^n \eta^p \eta^q C_{pq}^k \frac{\p}{\p \eta^k}, \sum_{a,b,c\in [n]} \phi^{abc}
\iota(e_a) \wedge \iota(e_b) \wedge \frac{\p}{\p \eta^c}\right]_\sS\\
&=& A + B
\Eeqrn
where
\Beqrn
A&=& \left[\sum_{l=1}^n \eta^l \iota(e_l),\sum_{a,b,c\in [n]} \phi^{abc}
\iota(e_a) \wedge \iota(e_b) \wedge \frac{\p}{\p \eta^c}\right]_\sS\\
&=& 2\sum_{a,b,c,l\in [n]} \eta^l\phi^{abc} [\iota(e_l), \iota(e_a)]\wedge \iota(e_b) \wedge \frac{\p}{\p \eta^c} +  \sum_{a,b,c\in [n]} \phi^{abc}
\iota(e_a) \wedge \iota(e_b) \wedge \iota(e_c)\\
&=& \sum_{a,b,c,l\in [n]} \eta^l(\phi^{kbc} C_{lk}^a +  \phi^{kca}C_{kl}^b)\iota(e_a)\wedge \iota(e_b)\wedge \frac{\p}{\p \eta^c}  +  \sum_{a,b,c\in [n]} \phi^{abc}
\iota(e_a) \wedge \iota(e_b) \wedge \iota(e_c)\\
\Eeqrn
and
\Beqrn
B&=& \left[\sum_{l=1}^n  - \frac{1}{2} \sum_{p,q,k=1}^n \eta^p \eta^q C_{pq}^k \frac{\p}{\p \eta^k}, \sum_{a,b,c\in [n]} \phi^{abc}
\iota(e_a) \wedge \iota(e_b) \wedge \frac{\p}{\p \eta^c}\right]_\sS\\
&=&\sum_{a,b,c,l\in [n]} \eta^l(\phi^{kca}C_{kl}^b)\iota(e_a)\wedge \iota(e_b)\wedge \frac{\p}{\p \eta^c} 
\Eeqrn
Thanks to the identity (\ref{2: inv of Cartan trivector}) many terms in $A+B$ cancel out and we conclude that
$$
-[\mu_1,\mu_3]= \iota(\phi) =\frac{1}{2}[\pi,\pi]_\sS=\frac{1}{2}[\mu_2,\mu_2]_\sS.
$$
\end{proof}

\subsubsection{\bf The first obstruction}
A straightforward calculation shows that $[\mu_3,\mu_3]_\sS\neq 0$ in general, more precisely
\Beq\label{2: [mu_3,mu_3]}
 [\mu_3,\mu_3]_\sS
 =4\sum_{\substack{a,b,c,i,j\\r,s}}
 \bigl(C_{rs}^{\ a}\phi^{rbi}\phi^{scj}\bigr)
 \iota_a\wedge \iota_b\wedge \iota_c\wedge 
\frac{\p}{\p \eta^i} \wedge  \frac{\p}{\p \eta^j} .
\Eeq
where we used the abbreviation
$$
\iota_a:=\iota(e_a).
$$
The precise form for $[\mu_3,\mu_3]_\sS$ is not important for the next deliberations, what is important is that $[\mu_3,\mu_3]_\sS$ is a cycle in $\wedge^\bu T\cM$ with respect to the differential 
\Beq\label{2: delta_S}
\delta_\sS:=[\mu_1,\ ]_\sS,
\Eeq
\Beqrn
\delta_\sS[\mu_3,\mu_3]_\sS&=&[\mu_1,  [\mu_3,\mu_3]_\sS]_\sS= 2[\mu_1, \mu_3]_\sS,\mu_3]_\sS=
 - [[\mu_2, \mu_2]_\sS,\mu_3]_\sS\\
&=& -2[[\mu_2, \mu_3]_\sS,\mu_3]_\sS\\
&=&0,
\Eeqrn
and has the form 
$$
 [\mu_3,\mu_3]_\sS
 =\sum_{\substack{a,b,c,i,j\\r,s}}
 \Phi^{abc|ij}
 \iota_a\wedge \iota_b\wedge \iota_c\wedge 
\frac{\p}{\p \eta^i} \wedge  \frac{\p}{\p \eta^j} .
$$
for some {\it constants}\, $\Phi^{abc|ij}$.

\subsection{A Lie subalgebra of $(\wedge^\bu T\cM, [ \ ]_\sS)$} Consider a subspace 
 $\cT^{p,q}_M$  of $\wedge^\bu T\cM$  spanned by polyvector fields  of the form
\Beq\label{2:Phi special}
\Phi=\sum_{a_1,...,a_p, i_1,..., i_q} \Phi^{a_1\ldots a_p\mid i_1\ldots i_q}
 \iota_{a_1}\wedge ...\wedge \iota_{a_p}\wedge
 \frac{\p}{\p \eta^{i_1}}\wedge ... \wedge \frac{\p}{\p \eta^{i_q}} \ \ p+q\geq 1.
\Eeq
for some {\it constants}\, $\Phi^{a_1\ldots a_p\mid i_1\ldots i_q}\in \K$, and set
$$
\cT_M^{spec}:= \bigoplus_{p+q\geq 1} \cT_M^{p,q}.
$$
  Note that both
$$
\mu_3 = -\sum_{a,b,c\in [n]} \phi^{abc} 
\iota_a \wedge \iota_b \wedge \frac{\p}{\p \eta^c} \ \ \text{and}\ \ \  [\mu_3,\mu_3]_\sS
$$
belong to the subspace $\cT^{spec}_M$. The latter observation about $[\mu_3,\mu_3]_\sS$ is not accidental.

\subsubsection{\bf Lemma}\label{2: Lemma on cT^spec} {\it (i) The subspace $\cT^{spec}_M\subset \wedge^\bu T\cM$ is closed with respect to the Schouten Lie bracket with
$$
\bigl[\cT_M^{p,q}, \cT_M^{r,s}\bigr] \subset \cT_M^{p+r-1, q+s}.
$$
(ii) For any quasi-Poisson structure $\mu_2\equiv\pi$ on $M$ one has}
$$
[\mu_2, \Phi]_\sS=0.
$$

\sip
\begin{proof} As 
$$
[\iota_a, \iota_b]_\sS= \sum_c C_{ab}^c \iota_c,\ \ [\iota_a, \frac{\p}{\p \eta^c}]_\sS=0,
$$
we have, for any
\Beqrn
 \Phi &=&\sum A^{a_1\ldots a_{p}\mid i_1\ldots i_{q}}
 \iota_{a_1}\wedge ...\wedge \iota_{a_{p}}\wedge
 \frac{\p}{\p \eta^{i_1}}\wedge ...\wedge  \frac{\p}{\p \eta^{i_{q}}},\\
 \Psi &=&\sum B^{b_1\ldots b_r\mid j_1\ldots j_s}
 \iota_{b_1}\wedge ...\wedge \iota_{b_r}
 \wedge \frac{\p}{\p \eta^{j_1}}\wedge ...\wedge\frac{\p}{\p \eta^{j_s}},
\Eeqrn
in $\cT_M^{spec}$,
\Beqrn
 [\Phi,\Psi]_{\sS}&=&\sum_{a_\bu,i_\bu,b_\bu,j_\bu} \sum_{u=1}^{p}\sum_{v=1}^{r}\pm C_{a_ub_v}^{c}
 A^{a_1... a_u... a_p\mid \ i_1\ldots i_{q}}
 B^{b_1... b_v... b_r\mid \mid j_1\ldots j_s }\iota_c\wedge
 \\
 && 
 \iota_{a_1}\wedge ...\widehat{\iota_{a_u}} ...\wedge \iota_{a_p}\wedge \iota_{b_1}\wedge ...\widehat{\iota_{b_v}} ... \wedge \iota_{b_r}
 \wedge  \frac{\p}{\eta^{i_1}}\wedge ... \wedge\frac{\p}{\p \eta^{j_s}}
 \Eeqrn
 which belongs to $\cT^{spec}$; hence $\cT^{spec}$ is a Lie subalgebra of $(\wedge^\bu T\cM, [ \ ,\ ])$.
Here the hat means that the corresponding factor is omitted.

\sip

As $\mu_2$ satisfies, by the very definition of a quasi-Poisson structure, the equations
$$
\bigl[\mu_2, \frac{\p}{\p \eta^{c}}\bigr]_\sS=0, \ \ \  \bigl[\mu_2, \iota_a\bigr]_\sS=0\ \ \ \forall\, a\in I, 
$$ 
we have
\Beqrn
 [\mu_2,\Phi]_{\sS}&=& \sum_{a_\bu, i_\bu}\sum_{k=1}^p
 \pm \Phi^{a_1... a_k... a_p\mid i_1...i_q}
 [\mu_2,\iota_{a_k}]_{\sS}
 \wedge \iota_{a_1}\wedge ...\widehat{\iota_{a_u}}...
 \wedge \iota_{a_p}\wedge\frac{\p}{\p \eta^{i_1}}\wedge ...\wedge  \frac{\p}{\p \eta^{i_{q}}}\\
 &=&0.
\Eeqrn
The Lemma is proven.
\end{proof}

The differential $\delta_\sS$ (see (\ref{2: delta_S})) in $\wedge^\bu T\cM$ does not respect the Lie subalgebra $\cT^{spec}_M$. Hence we extend temporarily the subspace $\cT^{p,q_M}$ to a larger subspace 
$\cT^{p,q|m}_M\subset \wedge^\bu T\cM$ which is generated, by definition, by polyvector fields of the form
$$
\Psi=\sum_{a_\bu, b_\bu i_\bu} \Phi^{a_1\ldots a_p\mid i_1\ldots i_q}_{b_1...b_m}\eta^{b_1} ... \eta^{b_m}
 \iota_{a_1}\wedge ...\wedge \iota_{a_p}\wedge
 \frac{\p}{\p \eta^{i_1}}\wedge ... \wedge \frac{\p}{\p \eta^{i_q}}.
$$
and let us set
$$
\cT^{ext}_M:= \bigoplus_{p+q\geq 1 \atop m\geq 0} \cT^{p,q|m}_M
$$
It is obvious that $\cT^{ext}_M$ is closed both with respect to the Schouten Lie bracket and, moreover, it is preserved by the differential $\delta_\sS$, i.e.\ $$
(\cT^{ext}_M,[\ ,\ ]_\sS, \delta_\sS)
$$
 is a dg Lie algebra. It contains $\cT^{spec}_M$ as the subspace with $m=0$.

\sip
Consider a  degree $-1$ linear map
$$
\fh: \cT^{p,q|m}  \lon  \cT^{p-1,q+1|m}
$$
defined on an arbitrary
\Beq\label{2: Phi p,q,m}
\Phi=\sum_{a_\bu, i_\bu} \Phi^{a_1\ldots a_p\mid i_1\ldots i_q}_{b_1...b_m}\eta^{b_1}... \eta^{b_m}
 \iota_{a_1}\wedge ...\wedge \iota_{a_p}\wedge
 \frac{\p}{\p \eta^{i_1}}\wedge ...\wedge \frac{\p}{\p \eta^{i_q}}  \in \cT^{p,q|m}
\Eeq
by the formula
\Beq\label{2: fh}
  \fh(\Phi):=\sum_{a_\bu, i_\bu}\sum_{k=1}^p(-1)^{m+k-1} \Phi^{a_1...a_k... a_p\mid i_1\ldots i_q}_{b_1...b_m}\eta^{b_1}... \eta^{b_m}
 \iota_{a_1}\wedge ...  \widehat{\iota_{a_k}} ... \wedge \iota_{a_p}\wedge
 \frac{\p}{\p \eta^{a_k}}\wedge \frac{\p}{\p \eta^{i_1}}\wedge ... \wedge \frac{\p}{\p \eta^{i_q}}  
\Eeq
It satisfies obviously the condition $\fh^2=0$. Moreover, $\fh(\cT^{spec})\subset \cT^{spec}$.

\subsection{An important change of viewpoint}\label{2: section on change of view} At this stage it is a good idea to forget about both the given $G$-manifold $M$ and the given representation
\Beq\label{2: iota from g to TM}
\Ba{rccc}
\iota: & \fg &\lon & TM\\
       & e_a &\lon & \iota_a:= \iota(e_a),
\Ea
\Eeq
and make our recursive construction of all the higher summands $\mu_{2n+1}$, $n\geq 1$, in $\mu_\hbar$ {\it universal}, i.e.\ independent of a particular choice of a $G$-manifold $M$ and a choice of a quasi-Poisson structure on $M$.

\sip

The spaces of polyvector fields $\cT^{spec}$, $\cT^{ext}$, their Lie bracket and the differential $\delta_\sS$ are independent of $M$ once we understand the symbols $\iota_a$ as some basis elements of the Lie algebra $\fg$, not as their representations as vector fields on some manifold. These spaces (and its elements $\mu_3$, $\mu_5$, etc which we construct soon) are mapped into $\wedge^\bu\cM$ under an arbitrary morphism of Lie algebras (\ref{2: iota from g to TM}).

\sip

From now on we treat the dg Lie algebra  $(\cT^{ext}, [\ , \ ]_\sS, \delta_\sS)$ and its subalgebra $\cT^{spec}$ purely formally as follows. Let $\fg$ be a Lie algebra with Lie bracket (\ref{2: Lie bracket in fg}), and consider a graded Lie algebra
$$
\fG:= \fg \oplus \fg[1] 
$$
whose basis, in the notation of \S {\ref{2: subsec about G and fg}}, is given by
$$
\{e_a, \frac{\p}{\p \eta^a} \}_{a\in I}
$$
so that $|e_a|=0, | \frac{\p}{\p \eta^a}|=-1$. 
The Lie brackets in $\fG$  are given in generators as follows,
$$
[e_a,e_b]=\sum_{c}C_{ab}^c e_c,\ \ \ \ [\frac{\p}{\p \eta^a}, \frac{\p}{\p \eta^b}]=[e_a, \frac{\p}{\p \eta^a}]=0.
$$
This Lie bracket extends naturally to the graded commutative algebra
$$
\odot^\bu (\fG[-1])
$$
 making the latter into a Gerstenhaber algebra.
Up to constants, this is precisely the formal Lie algebra $\cT^{spec}$ from the previous section,
\Beq\label{2: cT spec as wedge fg odot fg}
\cT^{spec}= \odot^{\bu\geq 1}(\fG[-1])\simeq  \bigoplus_{p+q\geq 1} \underbrace{\wedge^p\fg  \otimes \odot^q \fg[q]}_{\cT^{p,q}}.
\Eeq
In this incarnation as a Gerstenhaber algebra,  $\cT^{spec}$ can be identified with the subspace of the non-constant polynomials $\Phi(\iota_a, \psi_b)$ in the graded commutative algebra $\K[\iota_a, \psi_b]$. As $|e_a|=0$, $|\frac{\p}{\p \eta^a}|=-1$, the generators $\iota_a:=e_a[-1]$ and $\psi_a=\frac{\p}{\p \eta^a}[-1]$ of this graded commutative algebra have degees
$$
|\iota_a|=1,  \ \ \ \ |\psi_a|=0
$$, 
and the Schouten Lie bracket $[\ ,\ ]_S$  in $\cT^{spec}$ can be given now  by a Poisson type formula,
$$
\{ \Phi, \Psi\}_\sS= \sum_{a,b,c}  (-1)^{|\Phi|}
 \iota_c C^c_{ab}\frac{\p\Phi}{\p \iota_a} \frac{\p \Psi}{\p \iota_b}, \ \ \ \ \  \forall \Phi,\Psi\in \cT^{spec}.
$$
The formal version $\cT^{ext}$ of the larger Lie algebra $\cT^{ext}_M$
from the previous subsection can be identified with the following graded commutative polynomial algebra
$$
\cT^{ext}= \odot^\bu(\fg^*[-1])\ot \odot^{\bu\geq 1} (\fG[-1]) 
$$
which in turn can be identified with a subalgebra of the graded 
commutative algebra $\K[\iota_a, \psi_b, \eta^c]$, $|\eta^c=+1$,
spanned bu monomials of order $\geq 1$ in formal variables $\iota_a, \psi_b$.
The subspace $\cT^{p,q|m}$ is spanned by monomials of order $p$ in the formal variables $\iota_a$, of order $q$ in $\psi_b$ and of order $m$ in $\eta^c$. 
The Lie bracket $[\ ,\ ]_S$ in $\cT^{ext}_M$ also admit a nice Poisson-type presentation in $\cT^{ext}$  which is independent of $M$,
\Beq\label{2: Lie S  on T ext}
\{\Phi, \Psi\}_\sS=[\Phi, \Psi]_\sS= \sum_a 
\frac{\p\Phi}{\p \eta^a} \frac{\p \Psi}{\p \psi_a} +
 (-1)^{|\Phi|}\frac{\p\Phi}{\p \psi_a} \frac{\p \Psi}{\p \eta^a}+ \sum_{a,b,c}  (-1)^{|\Phi|}
 \iota_c C^c_{ab}\frac{\p\Phi}{\p \iota_a} \frac{\p \Psi}{\p \iota_b}
\Eeq
In this incarnation the  differential $\delta_\sS$ in $\cT^{ext}$ is given in terms of the Poisson bracket as follows,
\Beqr
\delta_\sS \Phi &=& \Biggl\{\underbrace{\sum_a \eta^a \iota_a - \frac{1}{2} \sum_{a,b,c} \eta^a \eta^b C_{ab}^c\psi_c}_{=:\nu_1}\ ,\  \ \ \  \Phi \Biggr\}_\sS \nonumber\\
&=& \sum_a \iota_a \frac{\p \Phi}{\p \psi_a}-
\sum_{a,b,c} \left(\eta^a \iota_c C_{ab}^c \frac{\p \Phi}{\p \iota_b} + \eta^b\psi_c C_{ab}^c\frac{\p \Phi}{\p \psi_a} - \frac{1}{2} \eta^a \eta^b C_{ab}^c\frac{\p \Phi}{\p \eta^c}\right) \label{2: delta_s in terms of S-brackets}
\Eeqr
while the operation $\fh$ defined in (\ref{2: fh}) becomes a simple derivation
\Beq\label{2: formal fh}
\fh(\Phi)= \sum_a \psi_a \frac{\p \Phi}{\p \iota_a}
\Eeq
After all these reformulations, a proof of the following crucial Lemma becomes a very short calculations.

\subsubsection{\bf Lemma}\label{2: Lemma on fh} {\it For any  $\Phi\in \bigoplus_{m\geq 0}\cT^{p,q|m}\subset \cT^{ext}$ one has
\Beq\label{2: (dh+hd)}
(\delta_\sS \circ \fh + \fh \circ \delta_\sS) \Phi = (p+q) \Phi.
\Eeq
Thus the complex $(\cT^{ext}, \delta_\sS)$ is acyclic.}

\begin{proof} 
Cancelling out second-order partial derivatives of $\Phi$ immediately, we obtain
\Beqrn
(\delta_\sS \circ \fh + \fh \circ \delta_\sS) \Phi &=&  \sum_{a,b,c} \eta^a \psi_c C_{ab}^c \frac{\p \Phi}{\p \iota_b} + \sum_a \iota_a \frac{\p \Phi}{\p \iota_a}
- \sum_{a,b,c} \psi_c\eta^a  C_{ab}^c \frac{\p \Phi}{\p \iota_b}+ \sum_a \psi_a \frac{\p \Phi}{\p \psi_a}
\\
&=&  \sum_a \iota_a \frac{\p \Phi}{\p \iota_a}  +\sum_a \psi_a \frac{\p \Phi}{\p \psi_a}= (p+q) \Phi
\Eeqrn
The Lemma is proven.
\end{proof}

\subsubsection{\bf Corollary}\label{2: corollary on delta_S exactness} {\it Every $\delta_\sS$-cycle $\Phi$ in $\cT^{p,q}\subset \cT^{spec}$ is $\delta_\sS$-exact, $\Phi=\frac{1}{p+q}\delta_\sS(\fh(\Phi))$.}

\sip

Let us emphasize again that the any representation (\ref{2: iota from g to TM}) induces the morphism of dg Lie algebras
$$
\iota: \bigl(\cT^{ext},\delta_\sS, \{\ ,\ \}_\sS\bigr)\lon
\bigl(\cT_M^{ext},\delta_\sS, [\ ,\ ]_\sS\bigr)
$$
and sends the map $\fh:\cT^{ext}\rar \cT^{ext} $ as defined in (\ref{2: formal fh}) to the map $\fh:\cT^{ext}_M\rar \cT^{ext}_M $ as defined in (\ref{2: fh}).

\subsection{Resolving the first obstruction}
If we set
\Beq\label{2: nu_3 in terms of phi}
\nu_3:= -\sum_{a,b,c\in [n]} \phi^{abc} 
\iota_a\iota_b \psi_c
\Eeq
then under the representation (\ref{2: iota from g to TM}) we have
$$
\mu_3=\iota(\nu_3)\ \ \text{and}\ \ \ [\mu_3,\mu_3]_\sS=\iota(\{\nu_3,\nu_3\}_\sS)
$$
If we define 
$$
\nu_5:=-\frac{1}{10}\fh\bigl([\nu_3,\nu_3]_\sS\bigr)=
-\frac65\sum_{
 \mathrm{Alt}_{abc}}
 \left(
t^{b\ell}t^{cv}
 C_{rs}^{a} t^{rk}C_{k\ell}^{i}t^{su}C_{uv}^{j}
 \right)
 \iota_b\wedge \iota_c\wedge \frac{\p}{\p \eta^a}\wedge \frac{\p}{\p \eta^i}\wedge\frac{\p}{\p \eta^j}
\in \cT^{2,3}
$$
then by Corollary {\ref{2: corollary on delta_S exactness}}
$$
\delta_\sS \nu_5 +\frac{1}{2}\{\nu_3, \nu_3\}_\sS=0.
$$
Hence setting
$$
\mu_5=\iota(\nu_5)=-\frac{1}{10}\fh\bigl([\nu_3,\nu_3]_\sS\bigr)
$$
we first conclude by Lemma {\ref{2: Lemma on cT^spec}}(ii) that
$$
[\mu_2,\mu_5]=0,
$$
and by the above,
$$
[\mu_1,\mu_5] + \frac{1}{2}[\mu_3,\mu_3]=0.
$$
Hence the ansatz
$$
\mu_\hbar^{<6}:= \mu_1 +\hbar\mu_2+ \hbar^2\mu_3 + \hbar^4\mu_5
$$
satisfies the equation
$$
[\mu_\hbar^{<6}, \mu_\hbar^{<6}]_\sS=0 \bmod O(\hbar^6).
$$
Now we are ready to beat all obstructions to constructing an MC element $\mu_\hbar$ in $(\wedge^\bu T\cM, [ \ \,\ ]_\sS)$.

\subsection{Proposition}\label{2: prop on collection nu_2n+1} {\it For any Lie algebra $\fg$ with an invariant scalar product as in \S {\ref{2: subsec about G and fg}} there exists a collection  of elements (called ``higher Jacobiators") 
\Beq\label{2: the collection nu_2n+1 n>1}
\{\nu_{2n+1}\in \cT^{spec}\}_{n\geq 1}
\Eeq
such that 
$$
\nu_3:= -\sum_{a,b,c} \phi^{abc} \iota_a\wedge \iota_b\wedge \p_{\eta^c}
$$
and for any $n\geq 2$ one has}
$$
\delta_\sS\nu_{2n+1}+ \frac{1}{2} \sum_{i+j= 2n+2\atop i,j\geq 3} 
\{\nu_i,\nu_j\}_\sS =0
$$
\begin{proof}
Starting with $\nu_3$ as above,
we assume by induction that for some $n\geq 2$ there is a sequence of elements,
$$
\left\{ \nu_{2k+1}\in \cT^{2, 2k-1}\right\}_{1\leq k\leq n}
$$
such that for each $k\in \{1,2,\ldots,n\}$ one has
$$
\delta_\sS\nu_{2k+1}+ \frac{1}{2} \sum_{i+j= 2k+2\atop i,j\geq 3} 
\{\nu_i,\nu_j\}_\sS =0.
$$
Define 
$$
\nu_{2n+3}:= -\frac{1}{2(2n+3)}\fh\left(\sum_{i+j= 2n+4\atop i,j\geq 3} 
\{\nu_i,\nu_j\}_\sS   \right).
$$
 As 
 \Beqrn
\delta_\sS \left(\sum_{i+j= 2n+4\atop i,j\geq 3}\{\nu_i,\nu_j\}_\sS  \right) 
&=&\sum_{i+j= 2n+4\atop i,j\geq 3}\bigl( \{\delta_\sS\nu_i,\nu_j\}_\sS 
- \{\nu_i,\delta_\sS\nu_j\}_\sS\bigr)\\
&=&-\sum_{a+b+k= 2n+5\atop a,b,k\geq 3} \{\{\nu_a,\nu_b\}_\sS, \nu_k\}_\sS \\
&=&0,
\Eeqrn
 we conclude by Corollary {\ref{2: corollary on delta_S exactness}}
 that 
 $$
\delta_\sS\nu_{2n+3}+ \frac{1}{2} \sum_{i+j= 2n+4\atop i,j\geq 3} 
\{\nu_i,\nu_j\}_\sS =0.
$$
 Hence defining recursively for any $n\geq 2$
 \Beq\label{2: recursion for nu_2n+1}
 \nu_{2n+1}:= -\frac{1}{2(2n+1)}\fh\left(\sum_{i+j= 2n+2\atop i,j\geq 3} 
[\nu_i,\nu_j]_\sS   \right).
\Eeq
we obtain a collection (\ref{2: the collection nu_2n+1 n>1}) with the desired properties. The proof is completed.
\end{proof}

\subsection{Invariance of every higher Jacobiator} By the very construction, every Jacobiator $\nu_{2n+1}$, $n\geq 1$, belongs 
to the subspace (cf. (\ref{2: cT spec as wedge fg odot fg})
$$
\nu_{2n+1}\subset \cT^{2,2n-1}\simeq \wedge^2\fg \ot \odot^{2n-1}\fg
$$
with the initial Jacobiator $\nu_3$ (which controls quasi-Poisson structures) belonging to $(\wedge^3)^\fg \subset \wedge^2\fg\ot \fg$.

\subsubsection{\bf Proposition} {\it Every Jacobiator is $\fg$-invariant}, i.e.\ for any $n\geq 1$ one has
$$
\nu_{2n+1}
 \in
 \left(\wedge^2\fg
 \ot \odot^{2n-1}\fg \right)^{\fg}.
$$
\begin{proof}   Applying the differential $\delta_\sS$ in $\cT^{ext}$ (see (\ref{2: delta_s in terms of S-brackets})) to an element $\Phi\in \cT^{spec}$ we obtain a general formula
\Beq\label{2: d_s applied to T spec}
 \delta_{\sS}\Phi
 \equiv \delta^{spec}\Phi -  \sum_a\eta^a\caL_a\Phi,
\Eeq
where the part
$$
\delta^{spec}:=\sum_a\iota_a\frac{\p}{\p\psi_a}
$$
is an endomorphism of $\cT^{spec}$, and the part
$-\sum_a\eta^a\caL_a$ with
$$
\caL_a:= \sum_{b,c}\left(
 \iota_c C^c_{ab}\frac{\p}{\iota_b}
 + \psi_c C^c_{ba}
 \frac{\p }{\p\psi_b}\right).
$$
sends $\Phi$ to the ambient space $\cT^{ext}\setminus \cT^{spec}$.
As each $\delta_S (\nu_{2n+1})$ belongs, by the very construction, to $\cT^{spec}$ we conclude that
$$
\caL_a (\nu_{2n+1})=0 \ \  \ \forall n\geq 1, \ a\in I.
$$
But the operator $\mathcal L_a$ is exactly the infinitesimal diagonal
$\fg$-action on
$$
 \wedge^p\fg
 \ot\odot^q\fg
 \simeq\cT^{p,q} \ \ \ \ \ \forall p+q\geq 1,
$$
Hence the claim follows.
\end{proof}

\subsection{Main Theorem}\label{2: Main Theorem on mu_hbar} {\it 
 For any quasi-Poisson manifold $(M, \pi, \iota)$, the associated graded manifold $\cM=M\times \fg[1]$ comes equipped with a graded Poisson structure
$$
\mu_\hbar=\mu_1 +  \hbar\mu_2 + \hbar^2\mu_3+ \sum_{n\geq 2} \hbar^{2n}\mu_{2n+1}\in \wedge^{\bu} T\cM[[\hbar]], \ \ \ 
\left[\mu_\hbar, \mu_\hbar\right]_{\sS}=0.
$$
where the first three terms are canonical,
\begin{align*}
\mu_1&= \sum_a\eta^a\iota_a
   -\frac12\sum_{a,b,c}(-1)^{|a|(|b|+1)}
       C_{ab}^{\ \ c}\eta^a\eta^b\frac{\p}{\p \eta^c},
\\
\mu_2&=\pi
\\
\mu_3 &= -\sum_{a,b,c} \phi^{abc}
\iota(e_a)\wedge \iota(e_b)\wedge \frac{\p}{\p \eta^c}
\end{align*}
 and all the higher polyvector fields
$$
\mu_{2n+1}:=\iota(\nu_{2n+1})
$$
are defined as representations in   $\wedge^{\bu} T\cM[[\hbar]]$ of the  collection (\ref{2: the collection nu_2n+1 n>1}) constructed recursively in Proposition {\ref{2: prop on collection nu_2n+1}}.}

\begin{proof} 
It remains to notice that all such $\mu_{2n+1}$ with $n\geq 1$
 commute  with $\mu_2$ by Lemma {\ref{2: Lemma on cT^spec}}(ii). 
\end{proof}

\subsubsection{\bf Corollary}\label{2: Lie_infty on cM} {\it The poly-differential operators 
$$
\mu_1: \Fun(\cM) \rar \Fun(\cM)[1], \ \ \ \mu_2: \wedge_\R^2 \Fun(\cM) \rar \Fun(\cM),
$$
$$
\left\{\mu_{2n+1}: \wedge_\R^{2n+1} \Fun(\cM) \rar \Fun(\cM)[1-2n]\right\}_{n\geq 1}
$$
make $\Fun(\cM)$ into a $\Lie_\infty$-algebra. }

\mip

The standard multiplication on $\Fun(\cM)$ together with the above three polydifferential operations $(\mu_1,\mu_2,\mu_{2n+1})$ can be understood also as a $Poisson_\infty$-structure on 
$\Fun(\cM)$, that is, a representation of the minimal resolution $Poisson_\infty$ of the operad of Poisson algebras on $\Fun(\cM)$.

\subsection{Functorial properties}\label{2: subsec on functorial}  Let $(M,\pi_M)$ and $(N,\pi_N)$ be two $G$-manifolds equipped with quasi-Poisson structures, and let
$$
F: M\rar N
$$
be a $G$-equivariant diffeomorphism such that
$$
F_*(\pi_M)=\pi_N.
$$
It induces a diffeomorphism of the associated supermanifolds
$$
\cF: \cM_M=M\times \fg[1] \lon \cM_N=N\times \fg[1]
$$
which is the identity map on the factors $\fg[1]$.
Let 
$$
\mu^M_\hbar=\mu_1^M+\hbar\mu_2^M +\sum_{n\geq 1} \hbar^{2n}\mu_{2n+1}^M  \ \ \ \text{and} \ \ \ 
\mu^N_\hbar=\mu_1^N+\hbar\mu_2^N +   \sum_{n\geq 1} \hbar^{2n}\mu_{2n+1}^N
$$
 be the  dg Poisson structures on $\cM_M$ and $\cM_N$ associated to the given quasi-Poisson structures on $M$ and $N$ respectively (see \S {\ref{2: Main Theorem on mu_hbar}}).

\subsubsection{\bf Proposition} $\cF_*(\mu^M_\hbar)=\mu^N_\hbar$.
\begin{proof} Since $F:M\rar N$ is a $G$-equivariant diffeomorphism, we have a commutative diagram
$$
\begin{tikzcd}[row sep=large, column sep=large]
& \fg \arrow[dl, "\iota_M"'] \arrow[dr, "\iota_N"] & \\
TM \arrow[rr, "F_*"'] & & TN
\end{tikzcd}
$$
of morphisms of Lie algebras. Hence 
$$
F_*(\iota_M(e_a))= \iota_N(e_a),\ \ \ \ F^*(\eta^a)=\eta^a
$$
for any $a\in [n]$ implying the equalities
$$
\cF_*(\mu_1^M)=\mu_1^N, \ \ \ \cF_*(\mu_{2n+1}^M)=\mu_{2n+1}^N,
$$
because $\cF_*$ commutes with the homotopy $\fh$ and the Schouten bracket.

As $\cF_*(\mu_2^M)=\mu_2^N$ by the very assumption, the claim is proven.
\end{proof}

We conclude that {\it the association $(M,\pi)\rar (\cM, \mu_\hbar)$ stated in the Main Theorem {\ref{2: Main Theorem on mu_hbar}} is functorial}.

\subsection{On essential uniqueness of higher Jacobiators}\label{2: subsec on essential uniqueness} Proposition {\ref{2: prop on collection nu_2n+1}} proves gives us a recursive construction  of a family of higher Jacobiators (\ref{2: the collection nu_2n+1 n>1}).  In this subsection we study the ambiguity of that construction.

\sip

Recall that the lowest Jacobiator $\nu_3$ is assumed to be fixed once and forever by the Cartan-trivector (\ref{2: cartan trivector})-(\ref{2: phi_abc}).
Assume that the lower Jacobiators are fixed
\[
 \nu_3,\nu_5,\ldots,\nu_{2n-1}.
\]
The equation determining the next term is given by (see Proposition {\ref{2: prop on collection nu_2n+1}})
\Beq\label{2: d_s nu_2n+1 = ...}
\delta_\sS\nu_{2n+1}=- \frac{1}{2} \sum_{i+j= 2n+2\atop i,j\geq 3} 
\{\nu_i,\nu_j\}_\sS.
\Eeq
If $\nu_{2n+1}$ and $\nu'_{2n+1}$ are two solutions of this equation,  then
$$
 \chi_{2n+1}:=
 \nu'_{2n+1}-\nu_{2n+1}
 \in(\cT^{2,2n-1})^{\fg}\cap\ker\delta_{\sS}.
$$
Conversely, adding any element of this kernel to a fixed solution $\nu_{2n+1}$
produces another solution.  Thus the space $(\cT^{2,2n-1})^{\fg}\cap\ker\delta_{\sS}$ describes all possible
differences between solutions.

\sip

On the other hand, in $\nu_{2n+1}$ is a solution of (\ref{2: d_s nu_2n+1 = ...}), then for any 
$$
 \xi_{2n}\in(\cT^{1,2n})^{\fg}
$$
the associated substitution
$$
 \nu_{2n+1}
 \longmapsto
 \nu_{2n+1}+\delta_{\sS}\xi_{2n},
$$
gives us a ``non-essential'', or ``gauge-trivial'' change of $\nu_{2n+1}$.  Hence differences of solutions belonging
to the space
$$
 (\cT^{2,2n-1})^\fg
 \cap \delta_{\sS}\bigl((\cT^{1,2n})^{\fg}\bigr) = 
 \delta_{\sS}\bigl((\cT^{1,2n})^{\fg}\bigr)
$$
are gauge trivial. We conclude that gauge-nontrivial changes of the $(2n+1)$-th Jacobiator
$\nu_{2n+1}$ are controlled the following graded vector space 
$$
 H_{\mathrm{essential\ ambiguity}}^{2n+1}
 :=
 \frac{(\cT^{2,2n-1})^{\fg}\cap\ker\delta_{\sS}}{
 \delta_{\sS}\bigl((\cT^{1,2n})^{\fg}\bigr)}.
$$
 Therefore a nonzero class (if any) in $H_{\mathrm{essential\ ambiguity}}^{2n+1}$ would represent a genuinely
gauge-inequivalent choice of the higher Jacobiator.

\subsubsection{\bf Lemma}\label{2: Lemma on H_amb =0} {\it For any $n\geq 2$ one has}
$$
H_{\mathrm{essential\ ambiguity}}^{2n+1}=0.
$$
\begin{proof} Given an arbitrary $\chi_{2n+1} \in 
(\cT^{2,2n-1})^{\fg}\cap\ker\delta_{\sS}$, consider
$$
 \xi_{2n}
 := \frac{1}{2n+1}\fh(\chi_{2n+1})
 \in(\cT^{1,2n})^{\fg},
 $$
 where $\fh$ is the homotopy (\ref{2: formal fh}). As $\delta_\sS\chi_{2n+1}=0$, we can write by Corollary {\ref{2: corollary on delta_S exactness}},
 $$
 \chi_{2n+1}=\delta_{\sS}\xi_{2n}.
$$
Lemma is proven.
\end{proof}

Given any two families of the Jacobiators as Proposition {\ref{2: prop on collection nu_2n+1}},
$$
\{\nu_{2n+1}\}_{n\geq 1}, \ \ \ \{\nu'_{2n+1}\}_{n\geq 1}
$$
with $nu_3=\nu_3'$ fixed by (\ref{2: nu_3 in terms of phi}). 
Then Lemma {\ref{2: Lemma on H_amb =0}} implies that the associated formal series
$$
\nu_\hbar
 = \sum_{n\geq 1}\hbar^{2n}\nu_{2n+1},
 \ \ \ \ \  
 \nu'_\hbar
 = \sum_{n\geq 1}\hbar^{2n}\nu'_{2n+1},
$$ 
are related by a gauge transformation
$$
 \nu'_\hbar
 =
 e^{-\ad_{\Theta_\hbar}}\nu_\hbar
 +
 \frac{1-e^{-\ad_{\Theta_\hbar}}}
      {\ad_{\Theta_\hbar}}
 \bigl(\delta_{\sS}\Theta_\hbar\bigr),
$$
for some
$$
 \Theta_\hbar
 =
 \sum_{n\geq 2}\hbar^{2n}\theta_{2n},
 \ \ \ \ 
 \theta_{2n}\in
 \bigl(\cT^{1,2n}\bigr)^{\fg}.
$$

Essential uniqueness means uniqueness of the gauge-equivalence class
of the family of higher Jacobiators.

\subsection{\bf Examples of higher Jacobiators}\label{2: example sl_2} Consider the Lie algebra
$$
 \fg=\mathfrak{sl}_2(\K)
 =\mathrm{span}_{\K}\{h,e,f\}
$$
with brackets
$$
 [h,e]=2e,\ \ \ \ [h,f]=-2f,\ \ \ \ [e,f]=h.
$$
The invariant tensor (\ref{2: invariant tensor f}) has the form,
$$
 t=\frac12 h\ot h +e\otimes f + f\otimes e
$$
and the associated invariant scalar product satisfies
$$
 \langle h,h\rangle=2,
 \ \ \ \ \
 \langle e,f\rangle=\langle f,e\rangle=1.
$$
As in \S {\ref{2: section on change of view}}  we extend $\mathfrak{sl}_2(\K)$ to the graded Lie algebra
$$
\fG:= \underbrace{\mathfrak{sl}_2(\K)}_{\text{spanned by}\ h,\, f,\, e} \oplus  \underbrace{\mathfrak{sl}_2(\K)[1]} 
_{\text{spanned by}\ 
 \p_{\eta^h},\, \p_{\eta^f}, \, \p_{\eta^e}}
 $$
whose generators satisfy the relations as just above and
$$
[\p_{\eta^a}, \p_{\eta^b}]=0, \ \ [\p_{\eta^a},h]= [\p_{\eta^a},f]=[\p_{\eta^a},e]=0  \ \ \ \forall a,b\in \{h,e,f\}.
$$
The Lie algebra $\cT^{spec}$ is defined as the associated Gerstenhaber algebra
$$
\cT^{spec}:= \odot^{\bu\geq 1}(\fG[-1]).
$$
Hence $\cT^{spec}$ is the subspace of non-constant
polynomials in the graded commutative algebra
$$
\K[\iota_h,\iota_e, \iota_f, \psi_h,\psi_e, \psi_f]
$$ 
whose generators $\iota_h,\iota_e, \iota_f$ have degree $1$, and
generators $\psi_h,\psi_e, \psi_f$ have degree zero. The Lie brackets (\ref{2: Lie S  on T ext}) restricted to $\cT^{spec}$ are given by
\Beqrn
 \{\Phi,\Psi\}_{\sS}
 &=&(-1)^{|\Phi|}\left(
 2\iota_e \frac{\p\Phi}{\p\iota_h}
 \frac{\p \Psi}{\p \iota_e}
 -2\iota_e
 \frac{\p\Phi}{\p\iota_e}
 \frac{\p\Psi}{\p\iota_h}
 -2\iota_f
 \frac{\p\Phi}{\p\iota_h}
 \frac{\p\Psi}{\p\iota_f}\right.\\
 &&\left. +2\iota_f \frac{\p\Phi}{\p\iota_f}
 \frac{\p \Psi}{\p \iota_h}
 +\iota_h
 \frac{\p\Phi}{\p\iota_e}
 \frac{\p\Psi}{\p\iota_f}
 -\iota_h
 \frac{\p\Phi}{\p\iota_f}
 \frac{\p\Psi}{\p\iota_e}\right).
 \Eeqrn
The initial Jacobiator (\ref{2: nu_3 in terms of phi}) is given  by
$$
\nu_3= 2\left(\iota_h\iota_e\psi_f - \iota_h\iota_f\psi_e + \iota_e\iota_f\psi_h\right)
$$
so that
$$
\{\nu_3, \nu_3\}=8 \iota_h\iota_e\iota_f( \psi_h^2+4\psi_e\psi_f)
$$
As the homotopy operator (\ref{2: formal fh}) is just the sum of partial derivatives
$$
\fh= \psi_h \frac{\p}{\p \iota_h} + \psi_f \frac{\p}{\p \iota_f}
+ \psi_e \frac{\p}{\p \iota_e},
$$
one can  solve the recursion equation (\ref{2: recursion for nu_2n+1}) for $\nu_5$ by direct calculations (which were first  done by ChatGPT 5.6, checked by the author, and are omitted in this paper):
$$
\nu_5=-\frac{4}{5} (\iota_h\iota_e\psi_f - \iota_h\iota_f\psi_e + \iota_e\iota_f\psi_h)( \psi_h^2+4\psi_e\psi_f)
$$
Further similar calculations  (first done by  ChatGPT 5.6 and then checked by the author) show that $\nu_7$ is given by
$$
\nu_7:=\frac{16}{35} (\iota_h\iota_e\psi_f - \iota_h\iota_f\psi_e + \iota_e\iota_f\psi_h)( \psi_h^2+4\psi_e\psi_f)^2.
$$
After these results the author noticed an obvious pattern, and  then the induction argument assisted by calculations by  ChatGPT 5.6 produced eventually a general formula
$$
 \nu_{2n+1}
 =
 \frac{24^nB_{2n}}{(2n)!}
 \left(\iota_h\iota_e\psi_f - \iota_h\iota_f\psi_e + \iota_e\iota_f\psi_h
 \right)
 \left(\psi_h^2+4\psi_e\psi_f\right)^{n-1}
$$
for all the Jacobiators in the case $\fg=\mathfrak{sl}_2(\K)$.
We omit the details.

\subsection{On deformation quantization of quasi-Poisson manifolds}  The above association of a dg Poisson structure $\mu_\hbar$ on the graded manifold $\cM$ associated to a quasi-Poisson structure $\pi$ on a manifold $M$ opens a new approach to  what can be called the deformation quantization of a quasi-Poisson manifold $(M,\pi)$. Since the Kontsevich formality map \cite{Ko} applies  to graded Poisson structures on (super)manifolds equally well as to ordinary Poisson structures on manifolds \cite{CF, LSX}, we define a {\it deformation quantization}\, of a given quasi-Poisson manifold $(M,\pi)$ as the Kontsevich deformation quantization of the associated  graded Poisson structure $\mu$ on $\cM$. The latter makes $\Fun(\cM)$ into, in general, $A_\infty$-algebra. Existing quantizations of quasi-Poisson manifolds due to Enriquez–Etingof \cite{EE}, Halbout , Li-Bland–Ševera \cite{LS}  treat important but restricted classes of quasi-Poisson manifolds by producing associative algebra objects in a Drinfeld braided category. 
 We hope to develop this approach to the quantization of quasi-Poisson structures  elsewhere.

\bip

\bip

{\Large
\section{\bf A $Poisson_\infty$ extension of the Fock-Rosly construction}
}

\mip

\subsection{Moduli space of flat $G$-connections}
Let $\Sigma$ be a connected oriented genus $g$ surface with boundary $\p\Sigma$ having $m\geq 1$ connected components.

\sip

Let $G$ be a connected Lie group with Lie  algebra  $\fg$ equipped with an invariant scalar product $t$ (as in \S {\ref{2: subsec about G and fg}}). 
A flat connection on a principal $G$-bundle over $\Sigma$ determines -- under a fixed trivialization of that bundle at some point $x_0\in \Sigma$ --  a holonomy homomorphism
\[
\rho:\pi_1(\Sigma,x_0)\longrightarrow G.
\]
Changing the trivialization at the base point conjugates the representation.  Hence the moduli space of flat $G$-connections  is the quotient space
\[
\fM_\Sigma(G)=\Hom(\pi_1(\Sigma),G)/G.
\]
with respect to the adjoint action of $G$ on itself. This is often a singular space so in practice one prefers working with the algebra of $G$-invariant functions on the so called framed moduli space
$$
\Hom(\pi_1(\Sigma),G)\simeq G^{2g+m-1}
$$
which is a smooth manifold. The latter isomorphism holds true because
$\pi_1(\Sigma,x_0)$ is isomorphic to the free group $F_r$ on $r=2g+m-1$ generators.

\subsection{A useful extension of $\pi_1(\Sigma)$}\label{2: subsec on Pi_1(Sigma,V)}
There is a very useful enlargement of the space $\Hom(\pi_1(\Sigma),G)$
associated with a choice of a nonempty finite set
\[
V=\{v_1,\dots,v_m\}\subset\p\Sigma.
\]
such that every boundary component of $\Sigma$ contains at least one point
of this set.
Instead of studying homotopy loops at a single base point $x_0\in \p\Sigma$, it is nowadays a common practice to work with the fundamental groupoid
\[
\Pi_1(\Sigma,V)
\]
whose objects are the points of $V$ and its morphisms $\ga: v_i\to v_j$ are homotopy classes, relative the endpoints,
of paths $\ga$ from $v_i$ to $v_j$; the point $v_i$ (resp. $v_j$) is called the {\it source} (resp.\ {\it target}) of $\ga$ and is denoted often by $s(\ga)$ (resp., $t(\ga)$). Given two morphisms $\al: v_i\rar v_j$ and $\be: v_j\rar v_k$ in this groupoid, we denote by $\be\cdot \al$ their composition giving us a morphism $v_i\rar v_k$. Thus $s(\be\circ \al)=s(\al)$ and $t(\be\circ \al)=t(\be)$. 

\sip

A representation of the groupoid $\Pi_1(\Sigma,V)$ in \(G\) is a functor
\Beq\label{2: repr of the groupoid Pi_1}
\rho:\Pi_1(\Sigma,V)\longrightarrow G,
\Eeq
where $G$ is regarded as the one-object groupoid associated to the group
$G$.

\subsection{The framed moduli space $M_{\Sigma,V}(G)$}
A flat connection $A$ in a principal $G$-bundle over $\Sigma$ with a chosen trivialization at the points of \(V\)  gives a representation
$$
\Ba{rccc}
\rho: & \Pi_1(\Sigma,V) & \lon & G\\
& \ga &\longmapsto& \Hol_\ga.
\Ea
$$
Note that 
$$
\Hol_{\gamma_1\gamma_2}
=
\Hol_{\gamma_1}\Hol_{\gamma_2}.
$$
Conversely, any representation (\ref{2: repr of the groupoid Pi_1})
determines the corresponding flat bundle with its trivializations at $V$.
The space of all representations,
\[
M_{\Sigma,V}(G)
=
\Hom\bigl(\Pi_1(\Sigma,V),G\bigr).
\]
is a smooth manifold parameterizing pairs
$$
\left(\text{a flat principal $G$-bundle $p:P\rar\Sigma$ \ \text{and}\  a framing, i.e.\ a choice of one point $\hat{v}\in p^{-1}(v)\ \forall v\in V$}\right).
$$
It is called the {\it framed}\, moduli space of flat connections on $(\Sigma,V)$.
Note that every $\al\in \Mor\bigl(\Pi_1(\Sigma,V)\bigr)$ gives us a smooth map
\Beq\label{2: Hol_al map from M(Sigma,V) to G}
\Ba{rccc}
\Hol_\al:M_{\Sigma,V}(G) & \lon G\\
           \rho   &\lon &  \Hol_\al(\rho)=\rho(\al).
\Ea
\Eeq

\subsection{\bf The intrinsic $\fg^V$-action}
\label{2: subsec g_V_action and CE}
 There is an action of the Lie group 
 $$
 G^V=\prod_{v\in V}G^v=\{V \rar G\}
 $$
 on the manifold $M_{\Sigma,V}(G)$,
$$
\Ba{ccc}
G^V \times M_{\Sigma,V}(G) & \lon & M_{\Sigma,V}(G)\\
\bigl(g=(g_v)_{v\in V}, \rho \bigr)    &\lon &  g\cdot \rho , 
\Ea
$$
which is given on an arbitrary homotopy path $\ga\in \Mor\bigl(\Pi_1(\Sigma,V)\bigr)$ by
$$
 (g\cdot \rho)(\ga):= g_{t(\ga)}\rho(\gamma)g_{s(\ga)}^{-1}.
$$
One has an isomorphism of quotient spaces
$$
\fM_\Sigma(G)\simeq M_{\Sigma,V}(G)/G^V.
$$

\sip

The $G^V$-action induces a morphism of Lie algebras
$$
X:  \fg^V=\bigoplus_{v\in V}\fg^v \lon  TM_{\Sigma,V}(G)=\Der(\Fun_{M_{\Sigma,V}(G)})
$$
whose restriction to the summand $\fg^v=\fg$  denoted by $X^v$. 
For any path $\al\in\Pi_1(\Sigma,V)$ the differential of the associated holonomy map (\ref{2: Hol_al map from M(Sigma,V) to G}) gives us a linear map of tangent vector spaces
$$
(d_\rho \Hol_\al)_*: T_\rho M_{\Sigma,V}(G) \lon
T_{\rho(\al)}G.
$$
Then for any $\zeta=(\zeta_v)_{v\in V} \in \fg^V$ one has
\Beq\label{2: X(zeta) as X^L - X^R}
(d_\rho\Hol_\ga)\bigl(X(\zeta)|_\rho\bigr)
 =\bigl(\zeta_{t(\ga)}\bigr)^R_{\rho(\ga)}
  -\bigl(\zeta_{s(\ga)}\bigr)^L_{\rho(\ga)}
\Eeq
where, for an element $\zeta\in \fg$, we denote by $\zeta^R$ and $\zeta^L$ the associated right- and left invariant vector fields on $G$. This formula gives us a complete characterization of the vector field $X(\zeta)$ on the framed moduli space.

\subsection{\bf The Chevalley--Eilenberg complex}
The above $\fg^V$-action 
makes the algebra of smooth functions
$\Fun_{M_{\Sigma,V}(G)}$ is a $\fg^V$-module, and we can consider the associated  Chevalley--Eilenberg complex
$$
C^\bullet(\fg^V, \Fun_{M_{\Sigma,V}(G)})=\Hom_\R\bigl(\wedge^\bu \fg^V, \Fun_{M_{\Sigma,V}(G)}\bigr)\simeq 
\Fun_{M_{\Sigma,V}(G)}\otimes \Lambda^\bullet (\fg^V)^*.
$$
equipped with the differential $\mu_1$ 
as in \S {\ref{2: CE complex of G-manifolds}}.
To simplify calculations, it is useful to understand the graded commutative algebra $\cA_{\Sigma,V}=C^\infty(M)\otimes \Lambda^\bullet \fg^*$ as an algebra of smooth functions on the dg (super)manifold (cf. \S {\ref{2: CE as cM}})
$$
\cM_{\Sigma,V}(G):= M_{\Sigma,V}(G) \times \fg^V[1]
$$
where vectors of the factor $\fg^V[1]$ are assumed to sit in degree $-1$.

\sip

 Let $\{e_a\}_{a\in I}$ be a basis of $\fg$; we denote by  $\{e_a^{(v)}\}_{a\in I}$ the associated basis in the summand $\fg^{v}$ of $\fg^V$. The corresponding dual basis  $\{e^a_{(v)}\}$  of $(\fg^v)^*\subset (\fg^V)^*$ induces odd (degree +1) coordinates on the supermanifold  $\cM_{\Sigma,V}(G)$ which we denote from now on by $\{\eta^a_{(v)}\}_{a\in I}$,
$$
\eta^a_{(v)}\eta^b_{(w)}=-\eta^b_{(w)}\eta^a_{(v)}, \ \ \ \ \forall \ v,w\in V.
$$
If we denote (arbitrary) local coordinates on the even part $M_{\Sigma,V}(G)$ symbolically by $x$, then every smooth function on the supermanifold $\cM_{\Sigma,V}(G)$ is given by a finite sum over all possible products of odd coordinate,
$$
f(x,\eta)=\sum_{k\geq 0} \sum f_{a_1...a_k}(x)\eta^{a_1}_{(v_{p_1})} ... \eta^{a_k}_{(v_{p_k})}.
$$
Hence for any point $v\in V$ the odd left partial derivative $\p/\p\eta^a_{(v)}$ makes an obvious sense as an operator on the graded commutative algebra
$$\Fun_{\cM_{\Sigma,V}(G)}.
$$
In this notation the Chevalley--Eilenberg differential $\mu_1$ is given by the following odd vector field 
\Beq\label{3: mu_1 on Sigma_V}
\mu_1= \sum_{v\in V} \underbrace{ \left(\sum_{a\in I} \eta^a_{(v)} X(e_a^{(v)}) - \frac{1}{2} \sum_{a,b,c\in I} \eta^a_{(v)} \eta^b_{(v)} C_{ab}^c \frac{\p}{\p \eta^c_{(v)}}\right)}_{\mu_1^{(v)}}
\Eeq
on $\cM_{\Sigma,V}(G)$. As
$$
H^0(\Fun_{\cM_{\Sigma,V}(G)}, \mu_1))=\bigl(\Fun_{\cM_{\Sigma,V}(G)}\bigr)^{G^V},
$$
we can view the dg manifold $(\cM_{\Sigma,V}(G),\mu_1)$ as a derived model for the moduli space of flat $G$-connections on the surface $\Sigma$.

\subsection{A skeleton of $\Sigma$}\label{3: skeleton Ga of Sigma}
A \emph{skeleton} $\Gamma\subset\Sigma$ is an embedded ribbon graph $\Ga$ such that
\Bi
\item[(i)] the set of vertices $V(\Ga)$ of $\Ga$ is precisely  the set  $V\subset\partial\Sigma$ marked points;
\item[(ii)] if $E(\Ga)$ stands for set of edges of $\Ga$, then each  $\fe\in E(\Ga)$ is smoothly embedded away from its endpoints; moreover, it is assumed that each edge $e$ is equipped with a fixed direction so that it makes sense to talk about {\it source}\, $s(\fe)$  (resp. {\it target}\, $t(\fe)$); the case $s(\fe)=t(\fe)$ is not excluded; 
\item[(iii)] $\Sigma$ deformation retracts onto $\Ga$.
\Ei
The latter condition implies that edges of $\Ga$ freely generate the fundamental groupoid in the sense of the following diffeomorphism
$$
\Psi_\Gamma:  G^{\#  E(\Ga)}  \stackrel{\sim}{\lon}  M_{\Sigma,V}(G)
$$
which sends  a tuple $\{g_\fe\in G\}_{\fe\in E(\Ga)}$ to the representation with
$\Hol(\fe)=g_\fe$ for all $\fe\in E(\Ga)$. The gauge group 
$$
G^{V(\Gamma)}=\{g=(g_v)_{v\in V},   \}
$$  acts on such tuples $\{g_\fe\}_{\fe\in E(\Ga)}$ as follows:
$$
g\cdot \{g_\fe\}_{\fe\in E(\Ga)} :=
\{ g_{t(\fe)}
g_\fe
g_{s(\fe)}^{-1}\},
$$
i.e.\ each vertex acts simultaneously on every incident half-edge. Let $H(v)$ be the set of half-edges attached to a vertex $v\in V(\Ga)$. Since $v\in \p \Sigma$, every set $H(v)=\{\fh_1,\fh_2,\ldots\}$ is totally ordered with respect to the orientation of 
$\Sigma$ as in the following picture 

\begin{center}
\begin{tikzpicture}[scale=1]
\draw[thick] (-2,0)--(2,0);
\fill (0,0) circle (2.5pt);
\node[below] at (0,0) {$v$};
\foreach \a/\lab in {25/1,55/2,85/3,120/4}{
  \draw[thick] (0,0)--(\a:1.8);
  \node at (\a:2.05) {$\fh_{\lab}$};
}
\draw[->] (1.55,.15) arc[start angle=5,end angle=155,radius=1.55];
\end{tikzpicture}.
\end{center}

To every $X\in \fg$ and every half-edge $\fh\in H(v)$ of an edge $\fe$ attached to $v$ we associate a vector field $(X)_\fh$ on the copy $G^\fe$ of the group $G$ in the product $G^{\# E(\Ga)}$ defined as follows,
$$
(X)_\fh :=
\begin{cases}
X^R,&\text{if $\fh$ is an incoming half-edge to $v$},\\[2mm]
- X^L,&\text{if $\fh$ is an outgoing half-edge from $v$}.
\end{cases}
$$
where $X^R$ (resp.\ $X^L$) stands for the right-invariant (resp. left-invariant) vector field on $G^\fe$ associated to $X$. The Fock-Rosly quasi-Poisson  structure on $G^{\# E(\Ga)}$ is given explicitly by \cite{FR} (see also \cite{ANPS, LS, MT})
$$
\pi_{\FR}^\Ga= 
\frac12
\sum_{v\in V(\Ga)}
\sum_{\substack{\fh',\fh''\in H(v)\\ \fh'< \fh''}}
\sum_{a,b}
t^{ab}\,
(e_a)_{\fh'}\wedge (e_b)_{\fh''} \in \wedge^2 T(G^{\# E(\Ga)})
$$
where we used notation introduced in \S {\ref{2: subsec about G and fg}}.
Then Theorem {\ref{2: Main Theorem on mu_hbar}} associates to $\pi_{FR}^\Ga$ a dg Poisson 
structure 
$$
\mu^\Ga_\hbar:=\mu_1 +\hbar\mu^\Ga_2 +\hbar^2\mu_3 + \sum_{n\geq 2}\hbar^{2n} \mu_{2n+1}
$$
on the supermanifold
$$
G^{\# E(\Ga)}\times \fg^V[1]
$$
such that $\mu^\Ga_2:=\pi_{FR}^\Ga$, $\mu_1$ is given by (\ref{3: mu_1 on Sigma_V}), $\mu_3$ is given by the Jacobiator
$$
\mu_3=- \sum_{v\in V}\sum_{a,b,c} \phi^{abc} X(e_a^{(v)})\wedge X(e_b^{(v)}) \wedge \frac{\p }{\p \eta^c_{(v)}},
$$
 and the higher operations  depend only on the structure constants $C_{ab}^c$ and the scalar  metric $t^{ab}$, the so called  higher Jacobiators from \S {\ref{2: prop on collection nu_2n+1}-\ref{2: Main Theorem on mu_hbar}}.
 Note that $\mu_1$ and $\{\mu_{2n+1}\}_{n\geq 1}$ do not depend on the choice of a skeleton $\Ga$.

\sip

The  quasi-Poisson structure bivector 
$$
\pi_{\FR}:=(\Psi_\Ga)_*\pi_{\FR}^\Ga \in\wedge^2TM(\SV)
$$
on $M_\SV(G)$ does not depend on the choice of a skeleton $\Ga$; one has
$$
(\Psi_\Ga)_*\pi_{\FR}^\Ga = (\Psi_{\Ga'})_*\pi_{\FR}^{\Ga'}
$$
if $\Ga$ and $\Ga'$ are related by edges reversals and slides (see \S 2.3 in \cite{ANPS} for a detailed proof). The diffeomorphism 
$$
\Psi_\Ga: G^{\# E(\Ga)} \lon M_{\SV}(G)
$$
extends to a diffeomorphism of graded manifolds
$$
\Psi_\Ga: G^{\# E(\Ga)}\times \fg^V[1] \lon \cM_\SV(G)=M_\SV(G)\times \fg^V[1] 
$$
as identity on the odd factor. 
Then applying Theorem {\ref{2: Main Theorem on mu_hbar}} to $\pi_{\FR}^\Ga$ and using we obtain the following result.

\subsubsection{\bf Proposition} {\it The derived moduli space $\cM_\SV(G)$ of flat $G$-connections comes naturally equipped with a dg Poisson structure
$$
\mu_\hbar:=(\Psi_\Ga)_*(\mu^\Ga_\hbar)=Q +\hbar\pi_{\FR} +\hbar^2\mu_3 + \sum_{n\geq 2}\hbar^{2n} \mu_{2n+1}
$$
which does not depend on the choice of a skeleton $\Ga$. }

\sip

The dg Poisson structure $\mu$ on $\cM_\SV(G)$ can be deformation quantized by applying to $\mu$ the classical Kontsevich's formality map \cite{CF,Ko}. The result is a curved (in general) $A_\infty$-structure on the graded algebra of smooth functions $\Fun(\cM_\SV(G))$ which will be studied in more detail elsewhere.

\bip

\bip

{\Large
\section{\bf From quasi-BV to BV$_\infty^{com}$}
}

\mip

\subsection{\bf Lie superalgebras with odd product: coordinate conventions}\label{4: convetions on Lie superalgebras}  A vector superspace is a $\Z_2$-graded vector space over a field $\K$
$$
V=V_{\bar0}\oplus V_{\bar1}.
$$
For a homogeneous element $x\in V$ we denote its parity by $|x|\in \Z_2=\{\bar0,\bar1\}$. A linear map $f:V\rar W$ between superspaces is called even (resp., odd) if $f(V_{\bar0})\subset W_{\bar0}$ and 
$f(V_{\bar1})\subset W_{\bar1}$ (resp.,  $f(V_{\bar0})\subset W_{\bar1}$ and 
$f(V_{\bar1})\subset W_{\bar0}$). 

\sip

There is a parity reverse endofunctor $\Pi$ in the category of superspaces. Given a vector superspace $V=V_{\bar0}\oplus V_{\bar1}$, the parity-reversed superspace $\Pi V$ is, by definition, the same underlying vector space with the parities exchanged,
$$
(\Pi V)_{\bar0}=V_{\bar1},\qquad (\Pi V)_{\bar1}=V_{\bar0}.
$$
For a homogeneous vector $x\in V$,
$$
|\Pi x|=|x|+1\mod 2\Z.
$$
Thus an odd linear map $V\to W$ between vector superspaces may equivalently be viewed as an even map $\Pi V\to W$, or as an even map $V\to\Pi W$.

\subsubsection{\bf Lie superalgebras}
A Lie superalgebra $\fg$ is a vector superspace equipped with even bracket
$[\  ,\ ]: \ot^2 V \rar V$ satisfying the graded skewsymmetry condition,
\[
[x,y]=-(-1)^{|x||y|}[y,x],
\]
and the graded Jacobi identity.

\sip

 Fix a homogeneous basis $\{e_a\}_{a\in I}$ of $\fg$
and let $\{e^a\}$ be the homogeneous dual basis of $\fg^*$,
$e^a(e_b)=\delta^b_a$. We abbreviate from now $|e^a|=|e_a|=|a|\in \Z_2$.
In such a basis the structure constants of the Lie bracket are given by
$$
[e_a, e_b]=\sum_{c} C_{ab}^{\ \, c} e_c, \ \ \  C_{ab}^{\ \, c}=-(-1)^{|a||b|}C_{ba}^{\ \,c}\in \K.
$$
Note that $C^{\ \, c}_{ab}\neq 0$ only if $|c|=|a|+|b| \bmod 2\Z$.

\subsubsection{\bf Odd invariant scalar product}
An \emph{odd invariant pairing} is a non-degenerate odd map
$$
\langle-,-\rangle:\fg\otimes\fg\to\K
$$
which is graded symmetric
$$
\langle x,y\rangle=(-1)^{|x||y|}\langle y,x\rangle
$$
 and satisfies
\[
\langle[x,y],z\rangle
+
(-1)^{|x||y|}
\langle y,[x,z]\rangle
=0
\]
for any $x,y,z\in \fg$.
One defines its structure constants in  homogeneous basis $\{e_a\}$ of $\fg$ by
$$
t_{ab}=\langle e_a,e_b\rangle, \ \ \ \ \  t_{ab}=(-1)^{|a||b|}t_{ba}\in \K.
$$
As $t$ is odd, 
$$
 t_{ab}\neq 0 \ \ \text{only if}\ \ |a|+|b|=1 \mod 2\Z. 
$$
Hence $t_{ab}\neq 0$ only if  $|a||b|=\bar0 $,  
and the graded symmetry condition takes the form
$$
t_{ab}=t_{ba}.
$$
Let $(t^{ab})$ be the inverse matrix. One has
$$
t^{ab}=t^{ba} \ \  \Rightarrow \ \ (-1)^{|a|}t^{ab}=- (-1)^{|b|}t^{ba}. 
$$
Hence if we define
$$
T^{ab}:= (-1)^{|a|}t^{ab}
$$
then 
$$
T^{ab}= - T^{ba}=- (-1)^{|a||b|} T^{ba}.
$$
as $T^{ab}\neq 0$ only if $|a||b|=\bar0$.
We define the inverse to the odd scalar product tensor by
\[
t:=\sum_{a,b}T^{ab}e_a\wedge e_b= \sum_{a,b}(-1)^{|a|}t^{ab}e_a\wedge e_b\in(\wedge^2\fg)^\fg.
\]
The invariance of $t$ translates in coordinates into the equations
\Beq\label{4: invariance of T}
\sum_{p} \bigl(C^{\ \, a}_{cp}T^{pb} + (-1)^{|a||b|}T^{ap}C_{cp}^{\ \, b}\bigr)=0
\Eeq
for any fixed $a,b,c$.

\subsubsection{\bf The Cartan trivector}
The odd pairing gives us an {\it odd}\,  isomorphism 
$$
\Ba{rccc}
t^\#: &\fg^* &\lon &\fg\\
      & e^a &\lon & t^\#(e^a):= \sum_c t^{ac}e_c.
\Ea
$$
Using it one defines a linear map 
$$
\Ba{rccc}
\phi: & \ot^3 \fg^* & \lon & \K\\
   &  (x^*,y^*, z^*) & \lon & x^*([t^\#(y^*), t^\#(z^*)]).
\Ea
$$
It can be understood equivalently as an invariant tensor $\phi\in (\odot^3\fg)^\fg$.
In coordinates, one has
$$
\phi=\sum_{abc}\phi^{abc}e_a\odot e_b \odot e_c \in\Sym^3\fg.
$$
where
\Beqr\label{4: phi abd in term of T and C}
\nonumber \phi^{abc}&= &\sum_{p,q} t^{bp}t^{cq}C^{\ \, a}_{pq}= \sum_{p,q}  (-1)^{|b|} t^{ap}C^{\ \, b}_{pq} t^{qc}=
 \sum_{p,q}  (-1)^{|a|+|b|+|c|+1} T^{ap}C^{\ \, b}_{pq} T^{qc} \\
 &=& -\sum_{p,q} T^{ap}C^{\ \, b}_{pq} T^{qc}
\label{4: phi abd in term of T and C}
\Eeqr 
One has
$$
\phi^{abc}=(-1)^{|a||b|}\phi^{bac}= (-1)^{|b||c|}\phi^{acb}.
$$
The invariance of $\phi$ translates in coordinates into the equations
\Beq\label{4: inv of phi}
\sum_k\left(C_{lk}^{\ \ a}\phi^{kbc}
 +(-1)^{|l||a|}C_{lk}^{\ \ b}\phi^{akc}
 +(-1)^{|l|(|a|+|b|)}C_{lk}^{\ \ c}\phi^{abk}\right)=0.
\Eeq

\subsection{Definition of a quasi-BV algebra \cite{ANPS}}
Fix an odd-metric Lie superalgebra
\[
(\fg,\langle-,-\rangle)
\]
and its canonical tensor $
\phi\in(\Sym^3\fg)^\fg\subset U(\fg)$. An $\fg$-{\it quasi-BV algebra} consists of:
\Bi
\item[(i)] a unital supercommutative associative superalgebra $A$;
\item[(ii)] a morphism of Lie superalgebras
\Beq\label{4: iota from fg to Der(A)}
\iota:\fg\lon {\Der}(A);
\Eeq
whose natural extension to a morphism of graded associative algebras
$$
\iota_U: U(\fg) \lon \End_\K(A)
$$
is denoted by the same letter $\iota_U$;
\item[(iii)] an odd differential operator
$$
\Delta: A\lon A
$$
of order at most two 
such that $
\Delta(1)=0$ and $\Delta$ is $\fg$-invariant in the sense that
\Beq\label{4: inv of Delta}
[\Delta,\iota(e)]=0
\ \ \ \forall\, e\in\fg,
\Eeq
\item[(iv)] a quasi-BV identity
\Beq\label{4: quasi-BV identity}
\Delta^2=\iota_U(\phi)\equiv \sum_{a,b,c}\phi^{abc}
 \iota(e_a)\iota(e_b)\iota(e_c)
\Eeq
holds true.
\Ei
If \(\varphi=0\), then
$\Delta^2=0$
and the structure becomes an ordinary BV algebra.

\subsection{Definition of a $BV_\infty^{com}$ algebra}
A  $BV_{\infty}^{com}$-{\em algebra}\, is, by definition \cite{Kr}, a differential graded commutative algebra  $(V,\delta)$
equipped with a countable collections of homogeneous linear maps, $\{\Delta_a: V\rar V,\
    |\Delta_a|=1-2a\}_{a\geq 1}$, such that each operator $\Delta_a$ is of order $\leq a+1$ (with respect to the given multiplication) and the equations,
    \Beq\label{5: BV_comm equation for Delta}
    \sum_{a=0}^n \Delta_a \circ \Delta_{n-a}=0, \ \ \ \text{with}\ \Delta_0:=-\delta
    \Eeq
hold for any $n\in \N$. These equations are equivalent to one equation,
    $$
    \Delta_\hbar^2=0
    $$
for the formal power series of operators
$$
\Delta_\hbar:=\sum_{a=0}^\infty \hbar^a \Delta_a,
$$
where the formal power variable $\hbar$ is assigned degree $2$.
The case when all $\Delta_a$ except $\Delta_1$ vanish corresponds to a $BV$-algebra.

\sip

If $\cB\cV_\infty^{com}$ stands for the dg operad controlling $BV_{\infty}^{com}$-algebras, and $\cB\cV$ for the operad controlling $BV$-algebras, then there is a natural surjection
$$
\cB\cV_\infty^{com}\lon \cB\cV
$$
is proven in \cite{CMW} to be a quasi-isomorphism. Moreover,
there is a {\it one-to-one correspondence}\,  between $BV_\infty^{com}$-algebras  and strongly homotopy involutive Lie bialgebras (see \S 5.4 in \cite{CMW} and \S 2.4 in \cite{Me}).

\sip

In this section we work solely with $\Z_2$-graded version of $BV_\infty^{com}$-algebras, so all operators $\Delta_a$ are assumed to be odd, while the parameter $\hbar$ is assumed to be even.

\subsection{``Derived" model for a quasi-BV algebra}\label{4: subsec on Def of cA}
Given a quasi-BV algebra  $(A, \Delta, \phi)$ as above
 Property (ii) says that $A$ is a $\fg$ module so we can consider the associated Chevalley-Eilenberg complex
\Beq\label{4: def of cA}
\cA:=\cA\otimes \odot^\bu( \Pi\fg^*)
\Eeq
which is a dg commutative superalgebra equipped with the differential (cf.\ (\ref{2: mu_1}))
$$
\Delta_1= \sum_a\eta^a\iota_a
   -\frac12\sum_{a,b,c}(-1)^{|a|(|b|+1)}
       C_{ab}^{\ \ c}\eta^a\eta^b\frac{\p}{\p \eta^c},
$$
where  $\{\eta^a:=\Pi e^a\}$ is a homogeneous basis of $\Pi\fg^*$.  Define
\begin{align}
\Delta_2&:=\Delta  \label{2: Delta_2}
\\
\Delta_3 &:= -\sum_{a,b,c} \phi^{abc}
\iota(e_a)\iota(e_b)\frac{\p}{\p \eta^c}
\label{2: Delta_3}
\end{align}

{\subsection{Proposition}\label{4: prop on Delta hbar <4}{\it For any  quasi-BV algebra the associated odd operator of order $\leq 3$
\Beq\label{4: Delta<4}
\Delta_\hbar^{<4}:= \underbrace{\Delta_1}_Q + \hbar \underbrace{\Delta_2}_\Delta + \hbar^2 \Delta_3
\Eeq
on the associated CE superalgebra $\cA[[\hbar]]$ satisfies the equation}
$$
(\Delta_\hbar^{<4})^2=0 \bmod O(\hbar^4)
$$
\begin{proof}
All three operators $\Delta_1$, $\Delta_2$ and $\Delta_3$ are odd.  Therefore
\Beq\label{4: Delta hbar squared}
 \Delta_\hbar^2
 =\Delta_1^2
 +\hbar[\Delta_1,\Delta_2]
 +\hbar^2\bigl(\Delta_2^2+[\Delta_1,\Delta_3]\bigr)
 +\hbar^3[\Delta_2,\Delta_3]
 \bmod O(\hbar^4)
\Eeq
The relation $\Delta_1^2=0$ is true as $\Delta_1$ is the standard coordinate expression for the CE-differential (cf.\ (\ref{2: mu_1})). As $\Delta$ is $\fg$-equivariant we have,
$$
 [\Delta_1,\Delta_2]=0,
 \ \ \ \ 
 [\Delta_2,\Delta_3]=0.
$$
Next we compute 
$$
[\Delta_1,\Delta_3]= \sum_a[\eta^a\iota_a,\Delta_3] -
\frac12\sum_{a,b,c}(-1)^{|a|(|b|+1)}
       [C_{ab}^{\ \ c}\eta^a\eta^b\frac{\p}{\p \eta^c}, \Delta_3].
$$
 Let us rewrite $\Delta_3$ in the following form
\Beq\label{4: Delta_3}
\Delta_3=\sum_c D^c\frac{\p}{\p \eta^c} \ \ \ \ \ \text{where}\ \ \ \
D^c:=-\sum_{a,b} \phi^{abc}\iota(e_a)\iota(e_b),\ \ |D^c|=|c|.
\Eeq
and abbreviate $\iota_a:=\iota(e_a)$. Then we have 
\Beqrn
 [\eta^a\iota_a,\Delta_3] &=& \sum_{a,c}[\eta^a\iota_a, D^c\frac{\p}{\p \eta^c}] \\
&=&\sum_c D^c\iota_c + \sum_{a,c} \eta^a[\iota_a, D^c]\frac{\p}{\p \eta^c}\\
&=& -\iota_U(\phi) - \sum_{a,b,c} \eta^a  D^b C_{ba}^{\ \, c}\frac{\p}{\p \eta^c}.
\Eeqrn
where we used the fact that
\Beqrn
 [\iota_a,D^c] &=& \sum_{i,j}\phi^{apq}[\iota_a,\iota_i\iota_j]\\
 &=&\sum_{i,j,k} \phi^{ijc}C_{ai}^{\ \ k}\iota_k\iota_j
   +(-1)^{|a||i|}\phi^{ijc}C_{aj}^{\ \, k}\iota_i\iota_k\\
 &=& -\sum_b  D^b C_{ba}^{\ \, c}\ \ \text{after applying (\ref{2: inv of Cartan trivector}) and relabelling}.
\Eeqrn
Next we compute Lie bracket with the second half of $\Delta_1$,
\Beqrn
-
\frac12\sum_{a,b,c}(-1)^{|a|(|b|+1)}
       [C_{ab}^{\ \ c}\eta^a\eta^b\frac{\p}{\p \eta^c}, \Delta_3]
 &=& -
\frac12\sum_{a,b,c,d}(-1)^{|a|(|b|+1)}
       [C_{ab}^{\ \ c}\eta^a\eta^b\frac{\p}{\p \eta^c}, D^d\frac{\p}{\p \eta^d}]\\
&=& \sum_{a,b,c} \eta^a D^b C_{ba}^{\ \, c}\frac{\p}{\p \eta^c}\\
\Eeqrn 
Summing up the results if the two calculations we conclude that
\Beqrn
[\Delta_1,\Delta_3] &=& -\iota_U(\phi) - \sum_{a,b,c} \eta^a  O^b C_{ba}^{\ \, c}\frac{\p}{\p \eta^c} +  \sum_{a,b,c} \eta^a  O^b C_{ba}^{\ \, c}\frac{\p}{\p \eta^c}\\
&=&  -\iota_U(\phi)=-\Delta_2^2
\Eeqrn
The vanishing of the coefficient of $\hbar^2$ in (\ref{4: Delta hbar squared}) is proven.

\sip

The vanishing of the coefficient of $\hbar^3$ in (\ref{4: Delta hbar squared}) is equal to
\Beqrn
[\Delta_2,\Delta_3] &=& \sum_c [\Delta, D^c] \frac{\p}{\p \eta^c}\\
&=& -\sum_{a,b,c}\phi^{abc}\left([\Delta, \iota_a]\iota_b + (-1)^{|a|}\iota_a[\Delta, \iota_b]\right) \frac{\p}{\p \eta^c}\\
&=& 0
\Eeqrn
where we used the $\fg$-invariance (\ref{4: inv of Delta}) of $\Delta$.

\sip

As every coefficient in (\ref{4: Delta hbar squared}) now vanishes, the proof of the identity
$(\Delta_\hbar^{<4})^2=0 \bmod O(\hbar^4)$ is completed.
\end{proof}

\subsection{The first obstruction for $\Delta_\hbar^{<4}$ to represent a $BV_\infty^{com}$ structure}
In general $(\Delta_\hbar^{<4})^2\neq 0$ because
\Beq\label{4: Delta_3^2}
\Delta_3^2 = \frac12\sum_{c,d}(-1)^{(|c|+1)|d|}
 [D^c,D^d]\, \frac{\p^2}{\p \eta^c\p \eta^d}
\Eeq
and
\Beqr
[D^c,D^d] &=&   \frac12\sum_{a,b,i,j} \phi^{abc}\phi^{ijd}
 [\iota_a\iota_b, \iota_i\iota_j] \nonumber\\
&=&\frac12\sum_{a,b,i,j,k} \phi^{abc}\phi^{ijd}\bigl( C_{bi}^{\ \, k}\iota_a \iota_k \iota_j
 +(-1)^{|b||i|}C_{ai}^{\ \, k}\iota_k \iota_b\iota_j
 +(-1)^{(|a|+|b|)|i|}C_{bj}^{\ \, k}\iota_i\iota_a\iota_k \nonumber\\
 &&
 +(-1)^{(|a|+|b|)|i|+|b||j|}C_{aj}^{\ \, k}\iota_i \iota_k\iota_b\bigr) \label{4: Delta_3^2 II}
\Eeqr
which does not vanish in general.
However, this calculation shows that the first obstruction 
$\Delta_3^2$ belongs to a rather special class $\f^{spec}_\cA\subset \Hom_\K(\cA,\cA)$ of  differential operators on $\cA$. More precisely, we define  $\f^{p,q}_\cA$ as a linear span over $\K$ of differential operators of the form (cf.\ \S 2)
\Beq\label{4: A p,q}
  \Sym(\iota_{a_1}\cdots\iota_{a_p})
 \p_{\eta^{i_1}}\cdots\p_{\eta^{i_q}},\ \ \ \  p+q\geq 1,              
\Eeq
and  set
$$
\f^{spec}_\cA=\bigoplus_{p+q\geq 1} \f_\cA^{p,q}.
$$
This superspace is a superalgebra with respect to the standard composition of operators on  $\cA$. The associated commutator Lie bracket satisfies,
$$
[ \f^{p_1,q_1}_\cA, \f^{p_2,q_2}_\cA]\subset \bigoplus_{k=0}^{p_1+p_2-1}  \f^{k,q_1+q_2}_\cA
$$

\subsubsection{\bf Proposition}\label{4: prop on [Delta, Phi]=0}
{\it For any quasi-BV structure $(A,\Delta, \iota)$ and any $\Phi\in \f^{spec}_\cA$ one has}
$$
[\Delta, \Phi]=0.
$$
\begin{proof}
This follows from the fact that $\Delta$ does not depend on variable $\eta^a$ and the invariance condition (\ref{4: inv of Delta}).
\end{proof}

\sip

The first obstruction $\Delta_3^2$ belongs to
$\f^{spec}_\cA$. Moreover, $\Delta_3^2$ is a cycle with respect to the differential 
\Beq\label{4: delta}
\delta:=[\Delta_1,...]=\bigl[ \sum_a\eta^a\iota_a-\sum_{a,b,c}
 \frac12(-1)^{|a|(|b|+1)}C_{ab}^{c}\eta^a\eta^b\partial_{\eta^c}
, \ldots\bigr]
\Eeq
 on $\Hom_\K(\cA,\cA)$. Indeed
 \Beqr
 \delta \Delta_3^2 &=&  [\Delta_1,\Delta_3]\Delta_3  -\Delta_3[\Delta_1,\Delta_3]  \label{4: Delta_3^2 acyclicity}\\
 &=&(-\Delta_2^2)\Delta_3  +\Delta_3(\Delta_2^2) \nonumber \\
 &=&-\Delta_2[\Delta_2,\Delta_3] + [\Delta_3,\Delta_2]\Delta_2 \nonumber\\
 &=& 0. \nonumber
 \Eeqr
The space $\f^{spec}_\cA$ is not closed under the differential $\delta$ so we enlarge it by first defining  $\f^{p,q\mid m}_\cA$ to be the span over $\K$ of the operators
$$
 \eta^{b_1}\cdots\eta^{b_m}
 \operatorname{Sym}(\iota_{a_1}\cdots\iota_{a_p})
 \partial_{\eta^{i_1}}\cdots\partial_{\eta^{i_q}}.
$$
and then setting
\Beq \label{eq:Aext-positive}
 \f^{ext}_\cA:=\bigoplus_{p+q\geq1\atop m\geq0}\;
   \f^{p,q\mid m}_\cA.
\Eeq
It contains $\f^{spec}_\cA$ as the  subspace,
$$
 \f^{spec}_\cA= \bigoplus_{p+q\geq1}\f^{p,q\mid 0}_\cA \subset \f^{ext}_\cA.
$$
The superspace $\f^{ext}_\cA$ is closed with respect to both the differential $\delta$  and the standard composition of operators, i.e.\ it is a dg superalgebra. 



\subsection{Making the construction universal}\label{4: subsec on making universal} The above formulae show that the only place in the definitions of Lie algebras  $\f^{spec}_\cA \subset \f^{ext}_\cA$ where we see a dependence on the choice of a particular quasi-BV algebra $A$ (and hence on its canonical extension $\cA$, see (\ref{4: def of cA})) is at the interpretation of derivations $\iota_a$ as images of the basis elements $e_a$ under the given morphism of Lie algebras; moreover, all we need to know about these symbols $\iota_a$ to keep our above statements about the pair  $\f^{spec}_\cA \subset \f^{ext}_\cA$ true is that they satisfy the commutativity relations
$$
 [\iota_a, \iota_b]=\sum_{c} C_{ab}^{\ \, c} \iota_c
$$
Hence it is a good idea at this stage to forget about a choice of a specific pre-BV algebra $\cA$ together with a particular representation (\ref{4: iota from fg to Der(A)}) and continue in a different formal setting as follows.

\sip

Consider an auxiliary Lie superalgebra (cf.\  \S {\ref{2: section on change of view}})
$$
{\fL}:= \fg \oplus \Pi\fg \oplus \Pi\fg^* \oplus \K
$$
whose basis, in the notation of \S {\ref{4: convetions on Lie superalgebras}}, is given by
$$
\{\iota_a=e_a,\ \ \p_{\eta^a}\simeq \Pi e_a, \ \ \eta^a\simeq \Pi e^a, \ \ t\in \K  \}_{a\in I}
$$
The Lie brackets in ${\fL}$  are defined  by declaring $t$ is central and by
$$
[\iota_a,\iota_b]=\sum_{c}C_{ab}^c \iota_c,\ \ [\p_{\eta^a}, \eta^b]=\delta^b_a t, \ \ 
 [\iota_a, \p_{\eta^b}]=0, \ [\iota_a, \eta^b]=0,\ \  [\p_{\eta^a},  \p_{\eta^b}]=0,  [{\eta^a},  {\eta^b}]=0.
$$

If in \S {\ref{2: section on change of view}}  we studied the Gerstenhaber algebra $\odot^{\bu \geq 1}(\fG[1])$ associated with the auxiliary graded Lie algebra $\fG$, in this section
we are interested in the quotient  $U({\fL})/(t-1)$ of the
universal enveloping algebra of the above superalgebra. The PBW theorem asserts an isomorphism of superspaces (not of superalgebras)
$$
U({\fL})/(t-1)\simeq \K[\iota_a,\p_{\eta^b}, \eta^b].
$$
We use this isomorphism to define the linear subspaces
$$
\f^{p,q|m}\subset U({\fL})/(t-1)
$$
as spanned by elements of the form
$$
 \eta^{b_1}\cdots\eta^{b_m}
 \operatorname{Sym}(\iota_{a_1}\cdots\iota_{a_p})
 \partial_{\eta^{i_1}}\cdots\partial_{\eta^{i_q}}.
$$
 Then
$$
\f^{ext}:= \bigoplus_{p+q\geq 1, m\geq 0} \f^{p,q|m}
$$
and 
$$
\f^{spec}:= \bigoplus_{p+q\geq 1} \f^{p,q|0}
$$
are sub-superalgebras of $U({\fL})/(t-1)$. In fact $\f^{spec}$ is the augmentation ideal $\bar{U}(\fg\oplus \Pi\fg)$ in the universal enveloping algebra 
$$
U(\fg \oplus \Pi\fg)\simeq U(\fg)\ot \odot^\bu(\Pi \fg). 
$$
of the subalgebra $\fg \oplus \Pi\fg\subset \fL$.

\sip

The superalgebra 
$\f^{ext}$ comes equipped with a differential given by the following supercommutator in $U({\fL})/(t-1)$,
\Beq\label{4: delta in f ext}
\delta \Phi:= \Biggl[\Delta_1:= \sum_a\eta^a\iota_a-\sum_{a,b,c}
 \frac12(-1)^{|a|(|b|+1)}C_{ab}^{c}\eta^a\eta^b\p_{\eta^c}, \ \Phi\Biggr] \ \ \forall \Phi\in \f^{ext}.
\Eeq

Given any quasi-BV algebra $A$
equipped with representation (\ref{4: iota from fg to Der(A)}), there is  an associated surjective morphism of dg superalgebras (denoted by the same letter)
\Beq\label{4:  iota frof f_ext}
\iota: (\f^{ext},\delta) \lon (\f^{ext}_\cA,\delta)
\Eeq
whose restriction gives us a surjective morphism of (non-differential) superalgebras
\Beq\label{4: iota fron f_spec}
\iota: \f^{spec} \lon \f^{spec}_\cA.
\Eeq
Under these morphisms the formal symbols $\iota_a$ become images $\iota(e_a)$ of basis vectors $e_a\in \fg$. 
 
 \subsection{Theorem} $H^\bu(\f^{ext},\delta)=0$
 \begin{proof}
 Consider an exhaustive filtration 
 $$
 0=F_0\subset F_1\subset\cdots \subset F_r\subset F_{r+1}\subset \ldots
 $$
 of the complex $(\f^{ext},\delta)$ with
 $$
 F_r:=\bigoplus_{\substack{1\leq p+q\leq r\\m\geq 0}}
 \f^{p,q\mid m}. 
 $$
 To prove the Theorem it is enough to show that the associated graded 
 complex
 $$
 \left(\mathrm{gr} \f^{ext}:=\bigoplus_{r\geq1}F_r/F_{r-1} ,  \delta_0\right)
 $$
 is acyclic. Any element $\Phi$ of $F_r/F_{r-1}\subset \mathrm{gr} \f^{ext}$ can be presented as a sum
 $$
 \Phi= \sum_{p+q=r\atop m\geq 0} \Phi^{a_1...a_p, i_1...i_q}_{b_1...b_m} \eta^{b_1}...\eta^{b_m} \Sym(\iota_{a_1}\cdots\iota_{a_p})
 \p_{\eta^{i_1}}\cdots\p_{\eta^{i_q}},
 $$
 where $\iota_a$ can be understood as a supercommuting (with everything) formal variables.
 
 The induced differential $\delta_0$ is given by
 $$
 \delta_0= \delta' + \delta''+ \delta'''
 $$
 where, for any $\Phi\in F_r/F_{r-1}$,
 \Beqrn
 (\delta'+\delta'')\Phi:&=& \bigl[\sum_a \eta^a\iota_a, \Phi\bigr] \bmod F_{r-1}\\
 &=&\underbrace{\sum_a \iota_a\frac{\p}{\p\, \p_{\eta^a}}\Phi}_{\delta'\Phi} + \underbrace{\sum_{a,b,c} (\text{Koszul sign})\, \eta^aC_{ab}^{\ \, c} \iota_c\frac{\p}{\p \iota_b}\Phi}_{\delta''\Phi} 
 \Eeqrn
 The action of the operator $\delta'$ is given in full details  by
 $$
 \delta'\Phi=\sum_{p+q=r\atop m\geq 0} (\text{Koszul sign})\, q\,\Phi^{a_1...a_p, a_{p+1}...i_q}_{b_1...b_m} \eta^{b_1}...\eta^{b_m} \Sym(\iota_{a_1}\cdots\iota_{a_p}\iota_{a_{p+1}})
 \p_{\eta^{i_2}}\cdots\p_{\eta^{i_q}},
 $$
 the action of $\delta''$ being defined similarly.
 The operator $\delta'''$ acts as the commutator
 $$
 \delta'''\Phi:=\bigr[ - \frac12(-1)^{|a|(|b|+1)}C_{ab}^{c}\eta^a\eta^b\p_{\eta^c}, \ \Phi\bigr] \bmod F_{r-1}.
 $$
Consider next an odd operator
$$
\fh_0:=\sum_a \p_{\eta^{a}} \frac{\p}{\p {\iota_a}}
$$
which acts on $\mathrm{gr} \f^{ext}$ in full details as follows
\Beq\label{4: fh_0}
\fh_0(\Phi):= \sum_{p+q=r\atop m\geq 0}(\text{Koszul sign})\, p\, \Phi^{a_1... a_{p-1}a_p...i_q}_{b_1...b_m} \eta^{b_1}...\eta^{b_m} \Sym(\iota_{a_1}\cdots\iota_{a_{p-1}})\p_{\eta^{a_p}} 
 \p_{\eta^{i_1}}\cdots\p_{\eta^{i_q}},
\Eeq
(Notice that $\fh_0^2=0$.)  Then we have (cf.\ \S {\ref{2: Lemma on fh}})
 \Beqrn
 (\delta_0\circ \fh_0 + \fh_0\circ \delta_0)\Phi 
 &=&\sum_a\left(
 \iota_a\frac{\p}{\partial \iota_a}
 +\p_{\eta^{a}}\frac{\p}{\p\, \p_{\eta^{a}}}\right)\Phi \\
 &=& (p+q) \Phi.
 \Eeqrn
 Thus every cycle in $(\mathrm{gr}\f^{ext},\delta_0)$ is a coboundary
 proving acyclicity of that complex. This fact together with the standard argument about a collapsing to zero converging spectral sequence of a filtered complex implies the acyclicity 
 of the original complex  $(\f^{ext},\delta)$.
 \end{proof}

\subsection{Extending the homotopy $\fh_0$ from $\mathrm{gr}\f^{ext}$ to $\f^{ext}$ (after \cite{GLS})  } As a (non-differential) superspace,  $\f^{ext}$ decomposes into the direct sum 
$$
\f^{ext}=\oplus_{r\geq 1} W_r
$$
where
$$
W_r:=\bigoplus_{\substack{ p+q= r\\m\geq 0}}
 \f^{p,q\mid m}. 
$$
The differential $\delta$ decomposes then into the sum
\Beq\label{4: delta splitting}
\delta=\delta_0 + \tau, 
\Eeq
where 
$$
\displaystyle \delta_0(W_r)\subset W_r\ \ \text{and}\ \  \tau(W_r)\subset \oplus_{k=1}^{r-1} W_{k}.
$$
We have
\Beq\label{4: split delta equations}
\delta_0^2=0, \ \ \ \delta_0\circ \tau + \tau\circ \delta_0 +   \tau^2=0.
\Eeq
Normalizing the action of the operator $\fh_0$ defined in (\ref{4: fh_0}) on elements $\Phi$ of $W_r$ by
$$
\wh{\fh}_0\Phi:=\frac{1}{r}\fh_0 \Phi \in F_r, 
$$
we get the identity
$$
(\delta_0\circ \wh{\fh}_0   +\wh{\fh}_0\circ \delta_0)\Phi=\Phi \ \ \ \ 
\forall \Phi\in \f^{ext}.
$$
Consider a linear map 
$$
\wh{\fh}: \f^{ext} \lon \f^{ext}
$$
defined by
$$
 \wh{\fh}(\Phi):=\sum_{k\geq 0}(-1)^k \wh{\fh}_0\circ (\tau\circ \wh{\fh}_0)^k\Phi \ \ \ \forall \Phi\in  \f^{ext}.
$$
Note that the sum above is finite as $\tau$ strictly lowers the value of the integer parameter $r$ used in the decomposition of $\f^{ext}$. It can be re-written as follows
$$
\wh{\fh}=\wh{\fh}_0\circ (1+ \tau\circ \wh{\fh}_0)^{-1}
$$
As
\Beqrn
\delta\circ (1+  \tau\circ \wh{\fh}_0) &=&  (\delta_0 + \tau)\circ (1+  \tau\circ \wh{\fh}_0)\\
&=& \delta_0+ \tau + (\delta_0\circ \tau +\tau^2) \circ \wh{\fh}_0\\
&=& \delta_0+ \tau - \tau\circ\delta_0 \circ \wh{\fh}_0\\
&=& \delta_0 + \tau\circ \wh{\fh}_0\circ\delta_0 \\
&=& (1+  \tau\circ \wh{\fh}_0)\circ \delta_0.
\Eeqrn
we have an equality
$$
\delta_0= (1+  \tau\circ \wh{\fh}_0)^{-1} \circ \delta\circ (1+  \tau\circ \wh{\fh}_0)
$$
Hence
\Beqrn
1 &=& \delta_0\circ \wh{\fh}_0   +\wh{\fh}_0\circ \delta_0\\
&=&(1+  \tau\circ \wh{\fh}_0)\ \circ (\delta_0\circ \wh{\fh}_0   +\wh{\fh}_0\circ \delta_0) \circ  (1+  \tau\circ \wh{\fh}_0)^{-1}\\
&=& \delta\circ \wh{\fh} +  \wh{\fh} \circ \delta.
\Eeqrn

Hence we proved the following

\subsubsection{\bf Proposition}\label{4: prop on cycles in f_spec}  {\it Any cycle $\Phi\in \f^{spec}$ satisfies the equality
\Beq\label{4: cycle in f_spec}
\Phi= \delta ( \wh{\fh}(\Phi)).
\Eeq
}

\subsection{Overcoming the first obstruction on the way to $BV_\infty^{com}$} Let $\cA$ be the superalgebra $A\ot \odot (\Pi \fg^*)$ associated with the given quasi-BV algebra $(A, \Delta, \iota)$
as in \S {{\ref{4: subsec on Def of cA}}.

\sip

Our purpose in this subsection is to find an order $\leq 5$ operator $\Delta_5$ on  $\cA$  such that the ansatz (\ref{4: Delta<4}) extends to
\Beq\label{4: Delta hbar <6}
\Delta_\hbar^{<6}=\Delta_\hbar^{<4} +\hbar^4 \Delta_5=  \Delta_1 + \hbar \Delta_2 + \hbar^2 \Delta_3+\hbar^4 \Delta_5
\Eeq
such that
\Beq\label{4: Delta<6 squared}
(\Delta_\hbar^{<6})^2=0 \bmod O(\hbar^6).
\Eeq
Moreover, we want to have a universal formula for $\Delta_5$, the one which works for any quasi-BV algebra.

\sip

The formal computations in (\ref{4: Delta_3}), (\ref{4: Delta_3^2}), and (\ref{4: Delta_3^2 acyclicity}) show  that 
the element
$$
\la_3:= -\sum_{a,b,c} \phi^{abc} \iota_a\iota_b \p_{\eta^c}\in \f^{2,1}\subset \f^{spec}
$$
has the following properties,
$$
\Delta_3= \iota (\la_3),\ \  \Delta_3^2= \iota (\la_3^2), \ \ \delta (\la_3^2)=0. 
$$
As $\la_3^2$ is a cycle in $\f^{spec}$, Proposition {\ref{4: prop on cycles in f_spec}} implies
$$
\la_3^2= \delta(\wh{\fh}( \la_3^2)).
$$
Hence defining
$$
\la_5:= -\frac{1}{2}\wh{\fh}([\la_3,\la_3])=-\wh{\fh}(\la_3^2), \ \ \ \Delta_5:=\iota(\la_5)
$$
we conclude that 
$$
\delta\la_5 + \frac{1}{2}[\la_3,\la_3]=0 \ \ \ \Rightarrow \ \ \ [\Delta_1,\Delta_5] + \frac{1}{2} [\Delta_3,\Delta_3]=0.
$$
By Proposition {\ref{4: prop on [Delta, Phi]=0}}
$$
[\Delta, \Delta_5]=0,
$$
and we finally conclude that the sum (\ref{4: Delta hbar <6}) with $\Delta_5$ defined just above satisfies the required equation (\ref{4: Delta<6 squared}).

\subsection{Proposition}\label{4: prop on collection nu_2n+1} {\it For any Lie superalgebra $\fg$ with an invariant odd scalar product as in \S {\ref{4: convetions on Lie superalgebras}} there exists a collection  of elements (``higher $\Delta$-Jacobiators") 
\Beq\label{4: the collection nu_2n+1 n>1}
\{\la_{2n+1}\in \f^{spec}\}_{n\geq 1}
\Eeq
such that 
$$
\la_3:= -\sum_{a,b,c} \phi^{abc} \iota_a\iota_b \p_{\eta^c}\in \f^{spec}
$$
and for any $n\geq 2$ one has}
$$
\delta\la_{2n+1}+ \frac{1}{2} \sum_{i+j= 2n+2\atop i,j\geq 3} 
[\la_i,\la_j] =0
$$
\begin{proof}
Starting with $\la_3$, let us assume by induction that
 for some $n\geq 2$ there is a sequence of elements,
\Beq\label{4: family la_2n+1}
\left\{ \la_{2k+1}\in \f^{spec}\right\}_{1\leq k\leq n}
\Eeq
such that for any $k\in \{2,\ldots,n\}$ one has
$$
\delta\la_{2k+1}+ \frac{1}{2} \sum_{i+j= 2k+2\atop i,j\geq 3} 
[\la_i,\la_j] =0
$$
Define 
$$
\la_{2n+3}:= -\frac{1}{2}\wh{\fh}\left(\sum_{i+j= 2n+4\atop i,j\geq 3} 
[\la_i,\la_j]   \right).
$$
The sum inside the round brackets is a cycle belonging to $\f^{spec}$,
$$
\delta \left(\sum_{i+j= 2n+4\atop i,j\geq 3}[\la_i,\la_j]  \right) 
=-\sum_{a+b+k= 2n+5\atop a,b,k\geq 3} \bigl[[\la_a,\la_b], \la_k\bigr]=0.
$$
Hence by Proposition  {\ref{4: prop on cycles in f_spec}} we obtain
 $$
\delta\la_{2n+3}+ \frac{1}{2} \sum_{i+j= 2n+4\atop i,j\geq 3} 
[\la_i,\la_j] =0.
$$
We conclude that, defining recursively for any $n\geq 2$
 $$
 \la_{2n+1}:= -\frac{1}{2}\wh{\fh}\left(\sum_{i+j= 2n+2\atop i,j\geq 3} 
[\la_i,\la_j]   \right),
$$ 
we obtain a collection (\ref{4: the collection nu_2n+1 n>1}) with the desired properties. The proof is completed.
\end{proof}

\subsection{\bf On the $\fg$-invariance of the $\Delta$-Jacobiators}
The superspace $\f^{spec}$ is identical to the non-constant part (which we mark by the subscript $+$) of the superalgebra $U(\fg)\ot \odot^\bu(\Pi\fg)$ (see \S {\ref{4: subsec on making universal}}). The latter space admits a natural diagonal $\fg$-action, and we denote 
the corresponding space of $\fg$-invariants by $(U(\fg)\ot \odot^\bu(\Pi\fg))^\fg$, and set
$$
(\f^{spec})^\fg:= (U(\fg)\ot \odot^\bu(\Pi\fg))_+^\fg.
$$

\subsubsection{\bf Proposition} {\it (i) If $\Phi\in \f^{spec}$ is a $\delta$-cycle, then
$$
\Phi\in (\f^{spec})^\fg.
$$

(ii) The homotopy $\wh{\fh}$ satisfies}
$$
\wh{\fh}\bigl( (\f^{spec})^\fg\bigr) \subset (\f^{spec})^\fg.
$$
\begin{proof} Every element of $\f^{ext}$ admits a unique representation as a linear combination of operators in the  ``normal" form
$$
\eta^{b_1}...\eta^{b_m} \Sym(\iota_{a_1}\cdots\iota_{a_p})
 \p_{\eta^{i_1}}\cdots\p_{\eta^{i_q}},
$$
This induces a direct sum decomposition into two summands
$$
\f^{ext}=\f^{spec} \oplus \f^{ext, \eta} 
$$
where $\f^{ext, eta}$ is the part containing at least one formal variable $\eta^a$. Let
$$
p_0: \f^{ext}\rar \f^{spec}, \ \ \ p_\eta: \f^{ext} \rar  \f^{ext, \eta}
$$
be the corresponding projections. Hence the differential (\ref{4: delta in f ext}) applied to an arbitrary $\Phi\in \f^{spec}$ split into a direct sum
$$
\delta(\Phi)=[\Delta_1, \Phi]= \underbrace{p_0[\Delta_1, \Phi]}_{\Delta^{spec}(\Phi)} \ + \ \underbrace{p_\eta[\Delta_1, \Phi]}_{\sum_a \eta^a\caL_a (\Phi)}
$$
where we use the fact that $p_\eta[\Delta_1, \Phi]$ contains only linear terms in the variables $\eta^a$. It is easy to see that 
$$
 \caL_a(\iota_b)=\sum_c C^c_{ab}\iota_c,
 \ \ \ \ 
 \caL_a(\p_{\eta^{b}})=(-1)^{|a|}\sum_cC^c_{ab}\p_{\eta^{c}},
$$
which corresponds to the diagonal $\iota_a$-action on the generators.
If $\Phi\in \f^{spec}$ is a cycle, then 
$$
\caL_a (\Phi)=0 \ \ \ \forall a\in I,
$$
implying that $\Phi\in (\f^{spec})^\fg$. The Claim(i) is proven.

\sip

Next we notice that 
$$
[\caL_a, \delta_0]=0, \ \ [\caL_a, \tau]=0, \ \ \  [\caL_a, \wh{\fh}_0]=0,
$$
implying 
$$
[\caL_a, \wh{\fh}]=0  \ \ \ \ \forall a\in I.
$$
This observation proves Claim(ii) as
$$
(f^{spec})^\fg=\{\Phi\in \f^{spec} | \caL_a(\Phi)=0\ \forall a\in I\}.
$$
\end{proof}

\subsubsection{\bf Corollary} {\it Each $\Delta$-Jacobiator from the family (\ref{4: family la_2n+1}) is invariant, i.e.
$$
\la_{2n+1}\in (\f^{spec})^\fg \ \ \ \forall n\geq 1.
$$
}

\mip

Note that the operator  $\fh_0$ preserves the weight (or the order) $r=p+q$ of the operators, the operators $\tau$ lowers it, so that the homotopy $\wh{\fh}$ does not increase the weight. Hence the weight
of $\la_{2n+1}$ is less than or equal to $2n+1$ implying that the operator $\iota(\la_{2n+1})$ on $\cA$ is of order $\leq 2n+1$.

\subsection{Main Theorem}\label{4: Theorem on Delta_hbar} {\it 
For any quasi-BV superalgebra $(A, \Delta, \iota)$, the associated CE superalgebra $\cA[[\hbar]]$ 
comes equipped with a $BV_\infty^{com}$-structure
$$
\Delta_\hbar=\Delta_1 +  \hbar\Delta_2 + \hbar^2\Delta_3+ \sum_{n\geq 2} \hbar^{2n}\Delta_{2n+1}, \ \ \  \Delta_\hbar^2=0,
$$
where the first three terms are canonical,
\begin{align*}
\Delta_1&= \sum_a\eta^a\iota_a
   -\frac12\sum_{a,b,c}(-1)^{|a|(|b|+1)}
       C_{ab}^{\ \ c}\eta^a\eta^b\frac{\p}{\p \eta^c},
\\
\Delta_2&=\Delta 
\\
\Delta_3 &= -\sum_{a,b,c} \phi^{abc}
\iota(e_a)\iota(e_b)\frac{\p}{\p \eta^c}
\end{align*}
 and all the higher operations 
$$
\Delta_{2n+1}:=\iota(\la_{2n+1})
$$
are defined as representations in  $\End_\K(\cA)$ of a  collection (\ref{4: the collection nu_2n+1 n>1}) constructed recursively in Proposition {\ref{4: prop on collection nu_2n+1}}.}

\begin{proof} 
It remains to notice that, by Proposition {\ref{4: prop on [Delta, Phi]=0}}, all   the constructed operators
 $\Delta_{2n+1}$, $n\geq 1$,
 commute  with $\Delta_2$. 
\end{proof}

\subsection{Examples} Following \cite{ANPS} the above result applies to the framed moduli space $M_\SV(G)$ of flat $G$-connections on a surface $\Sigma$ with boundaries under the assumption that the Lie superalgebra $\fg$ of the given Lie supergroup $G$ has an odd invariant scalar product. More concretely, for a fixed ribbon graph skeleton $\Ga\subset \Sigma$  the authors of \cite{ANPS} constructed explicitly (see formula (10) ibid.)  a quasi-BV algebra structure on $\Fun(M_\SV(G))$ and then proved that $\Delta$ is independent of the choice of $\Ga$, i.e.\ it is intrinsic.  The Main Theorem {\ref{4: Theorem on Delta_hbar}} applies directly to the quasi-BV algebra $(\Fun(M_\SV(G)), \Delta)$ and hence gives an associated $BV_\infty^{com}$-structure $\Delta_\hbar$ on the associated supermanifold $M\times \Pi\fg$. Since the construction outlined in Main Theorem {\ref{4: Theorem on Delta_hbar}} is obviously functorial (cf.\ \S {\ref{2: subsec on functorial}}), the obtained $BV_\infty^{com}$-structure $\Delta_\hbar$ is also intrinsic.

\end{document}